\documentclass{article}
\usepackage{PRIMEarxiv}
\usepackage{amssymb}
\usepackage{amsmath}
\usepackage{flafter}
\usepackage{neuralnetwork}
\usepackage{afterpage}
\usepackage[english]{babel}
\usepackage{amsthm}
\newtheorem{prop}{Proposition}
\usepackage{hyperref}

\usepackage{subfigure}
\usepackage{graphics,graphicx}
\usepackage{enumerate}
\usepackage[figuresright]{rotating}
\usepackage[utf8]{inputenc} % allow utf-8 input
\usepackage[T1]{fontenc}    % use 8-bit T1 fonts
\usepackage{hyperref}       % hyperlinks
\usepackage{url}            % simple URL typesetting
\usepackage{booktabs}       % professional-quality tables
\usepackage{amsfonts}       % blackboard math symbols
\usepackage{nicefrac}       % compact symbols for 1/2, etc.
\usepackage{microtype}      % microtypography
\usepackage{lipsum}
\usepackage{fancyhdr}       % header
\usepackage{graphicx}       % graphics
\graphicspath{{media/}}     % organize your images and other figures under media/ folder

\title{Computing Anti-Derivatives using Deep Neural Networks
\thanks{\textit{\underline{Submitted to Elsevier}}: 
\textbf{}} 
}

\author{
  Dibyajyoti Chakraborty , Shivasubramanian Gopalakrishnan\\
 Department of Mechanical Engineering\\ Indian Institute of Technology, Bombay\\ Mumbai 400076\\ India\\
}

\begin{document}
\maketitle

\begin{abstract}
This paper presents a novel algorithm to obtain the closed-form anti-derivative of a function using Deep Neural Network architecture. In the past, mathematicians have developed several numerical techniques to approximate the values of definite
integrals, but primitives or indefinite integrals are often non-elementary. Anti-derivatives are necessarily required when there are several parameters in an integrand and the integral obtained is a function of those parameters. There is no theoretical method that can do this for any given function. Some existing ways to get around this are primarily based on either curve fitting or infinite series approximation of the integrand, which is then integrated theoretically. Curve fitting approximations are inaccurate for highly non-linear functions and require a different approach for every problem. On the other hand, the infinite series approach does not give a closed-form solution, and their truncated forms are often inaccurate. We claim that using a single method for all integrals, our algorithm can approximate anti-derivatives to any required accuracy. We have used this algorithm to obtain the anti-derivatives of several functions, including non-elementary and oscillatory integrals. This paper also shows the applications of our method to get the closed-form expressions of elliptic integrals, Fermi-Dirac integrals, and cumulative distribution functions and decrease the computation time of the Galerkin method for differential equations.
\end{abstract}

% keywords can be removed

\section{ Introduction} 
There has been some research to find methods\cite{muhammad2003double} that can approximate the anti-derivatives. Most of them are based on polynomial or exponential spline fitting, which often can not capture some integrals' highly non-linear and non-elementary nature. We can not use these methods directly in cases where the integrand is dependent on various parameters. They require a different approach for every integral. Neural networks, acting as universal approximators\cite{hornik1989multilayer}, can be a potent tool for this purpose. Definite integrals have been approximated \cite{zhe2006numerical} using a single hidden layer neural network to significant accuracy. Dual Neural Networks \cite{li2019dual} has also been used to calculate definite integrals in the cases where the integrand is represented as discrete values. Computational speedup over other numerical techniques using shallow neural networks has also been shown\cite{lloyd2020using}.
For the integral, 
\begin{equation}
    \int f(x)\,dx
\end{equation}
the previous works have shown the use of a single-layer neural network ( $N_1(x)=(w_j^T\sigma(w_i^Tx+b_i)+b_j)$) to approximate the integrand, similar to curve fitting. 
Later, this approximation of the integrand is integrated theoretically, which in a way beats the purpose. Since theoretically integrating deep neural networks is not feasible, they have used shallow neural networks that are often insufficient for highly non-linear functions.\\
\begin{equation}
    Older\ Methods\ :\ \int f(x)\,dx = \int N_1(x)\,dx
\end{equation}
This paper presents the algorithm- Deep Neural Network Integration(DNNI)
\begin{equation}
    DNNI : \int f(x,a,b,...)\,dx = N(x,a,b,...)
\end{equation} where $N(x,a,b,..)$ is a Deep Neural Network similar to figure \ref{fig 1:DNNI}. The breakthrough that has propelled the DNNI algorithm is automatic differentiation\cite{baydin2018automatic} which enables us to include the derivative of the neural network in its loss function. Using DNNI, we can obtain the anti-derivative directly as a continuous function without theoretically integrating.

DNNI can be a single method for approximating primitives, calculating the value of definite integrals, and obtaining the closed-form expression of an integral as a function of other parameters.\\ This paper is organized as follows: Section 2 describes the algorithm; In section 3, anti-derivatives have been computed using DNNI and compared with theoretical results for simple, complicated functions, non-elementary integrals, and oscillatory functions. Section 4 focuses on the applications of DNNI in obtaining the closed-form expressions approximating the elliptic and Fermi-Dirac integrals, finding cumulative distribution functions, and speeding up numerical differential equations solvers. The codes used in the paper are available at \href{https://github.com/Dibyajyoti-Chakraborty/Deep-Neural-Network-Integration}{https://github.com/Dibyajyoti-Chakraborty/Deep-Neural-Network-Integration}. They are executed in Intel® Core™ i7-9700K CPU @ 3.60GHz × 8 processors, 32 GB RAM and NVIDIA Corporation GV100GL [Tesla V100 PCIe 32GB] GPU.

\section{ Methodology}
Integration of a continuous real function f is defined as 
\begin{equation}
    I(x,a,b,..) = \int f(x,a,b,..)\,dx
\end{equation}
where $I$ is the anti-derivative of $f$ since 
\begin{equation}
   \frac{\partial I} {\partial x}=f(x,a,b,..)
\end{equation}

\begin{figure}[h]
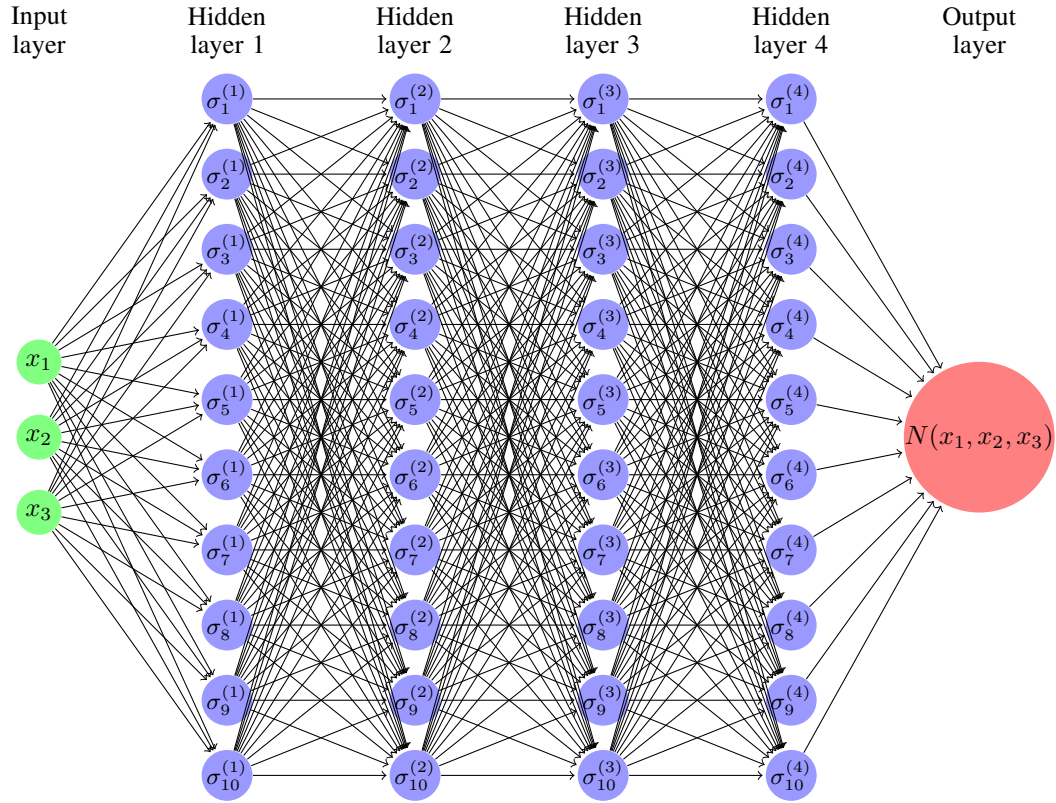

	\centering
\begin{neuralnetwork}[height=10]
        \newcommand{\x}[2]{$x_{#2}$}
        \newcommand{\y}[2]{$N(x_1,x_2,x_3)$}
        \newcommand{\hfirst}[2]{\small $\sigma^{(1)}_{#2}$}
        \newcommand{\hsecond}[2]{\small $\sigma^{(2)}_{#2}$}
        \newcommand{\hthird}[2]{\small $\sigma^{(3)}_{#2}$}
        \newcommand{\hfourth}[2]{\small $\sigma^{(4)}_{#2}$}
        \inputlayer[count=3, bias=false, title=Input\\layer, text=\x]
        \hiddenlayer[count=10, bias=false, title=Hidden\\layer 1, text=\hfirst] \linklayers
        \hiddenlayer[count=10, bias=false, title=Hidden\\layer 2, text=\hsecond] \linklayers
        \hiddenlayer[count=10, bias=false, title=Hidden\\layer 3, text=\hthird] \linklayers
        \hiddenlayer[count=10, bias=false, title=Hidden\\layer 4, text=\hfourth] \linklayers
        \outputlayer[count=1, title=Output\\layer, text=\y] \linklayers
    \end{neuralnetwork}
	\caption{Deep Neural Network architecture with three inputs and one output. It also has four hidden layers with ten nodes each.}
	\label{fig 1:DNNI}
\end{figure}

In the DNNI algorithm, the integral $I$ is approximated by a Feed Forward Deep Neural Network $N(x)$ as shown in figure \ref{fig 1:DNNI}. The Neural Network can be represented as
\begin{equation}
    N(x_1,x_2,...) = W_{L+1}^T\sigma(W_L^T(\sigma(W_{L-1}^T(\sigma(..........)+b_{L-1})+b_{L})+b_{L+1}
\end{equation} 
where $W$s and $b$s are weights and biases respectively and $\sigma$ is a non-linear function like sigmoid$(1/(1+e^{-x}))$ or tanh. Now, in a required domain $x\in [\alpha,\beta]$ which can be arbitrarily large, we aim to make
\begin{equation}
    N(x,a,b,...) \approx I(x,a,b,...) 
\end{equation}
\begin{equation}
    \implies N(x,a,b,...) \approx \int f(x,a,b,...)\,dx 
\end{equation}
\begin{equation}
    \implies \frac{\partial N(x,a,b,...)}{\partial x} \approx  f(x,a,b..) 
\end{equation}
Since a neural network is a continuous function, it can be easily differentiated using the latest developments in automatic differentiation. This differential of the neural network is included in the loss function to minimize its deviation from the integrand. As the derivative of the neural network gets closer to the integrand, the integral is approximated by the neural network. Hence, a loss function can be formed such as \begin{equation}
    LOSS :\ \  MSE\left(\frac{\partial N(x,a,b,...)}{\partial x} , f(x,a,b..)  \right)
\end{equation}
where
\begin{equation}
    MSE(x,y) :\ \ \frac{\sum_{i=1}^{N}(x_i-y_i)^2 }{N}
\end{equation}
is the Mean Squared Error. The weights and biases can be tuned using an optimization algorithm to reach the required accuracy. Gradient Descent or quasi-Newton based optimization algorithms are mostly used. In this paper, we have mostly used the Adam algorithm\cite{kingma2014adam}, with learning rate scheduling, for optimization. \\In the cases where the limits are defined, DNNI can give a closed form approximation of the integral as a function of other parameters. 
\begin{equation}
    \int_{x_0}^{x_n} f(x,a,b,...)\,dx = F(a,b,...) \approx N(x_n,a,b,...)-N(x_0,a,b,...)\end{equation}
The model depth, the number of nodes in each hidden layer, and the number of epochs are selected based on the complexity. Some activation functions used are shown in the table below. ReLU is not used as its higher gradients vanish.
\begin{center}
\begin{tabular}{ c c c c} 
  \hline  \hline
  Name &Expression & Derivative & Second derivative \\ 
  \hline  \hline
  Sigmoid    & $ \frac{1} {1 + e^{-z}}$ & $\frac{e^{x}} {(1+e^{x})^2}$ & $-\frac{e^{x}(e^{x}-1)} {(1+e^{x})^3}$ \\ 
  \hline
  Tanh & $tanh(x)$ & $sech^2(x)$ & $-2\ tanh(x)\ sech^2(x)$ \\ 
  \hline ReLU & max(0,x) & $
\begin{array}{cc}
  \Bigg\{  
    \begin{array}{cc}
      0 & x\leq 0 \\
      1 & 0\le x \\
    \end{array}
\end{array}
$ & 0\\\hline \hline
\end{tabular}
\end{center}

\section{ Results}
We tried the DNNI algorithm in many cases setting the lower limit to some arbitrary value to eliminate the constant of integration. Then, we compared the theoretical primitive, if available, with its DNNI approximation.
\subsection{Simple Integrals}
In the following subsection, we have shown the application of DNNI to obtain the anti-derivatives of some ubiquitous integrals. They have been plotted with their theoretical counterpart for comparison. \\
 Case 1:
\begin{equation}
    \int x^6 \,dx = x^7/7 + c
\end{equation}
\begin{figure}[h]
    \centering
    \includegraphics[width=7.5cm,height=5.6cm]{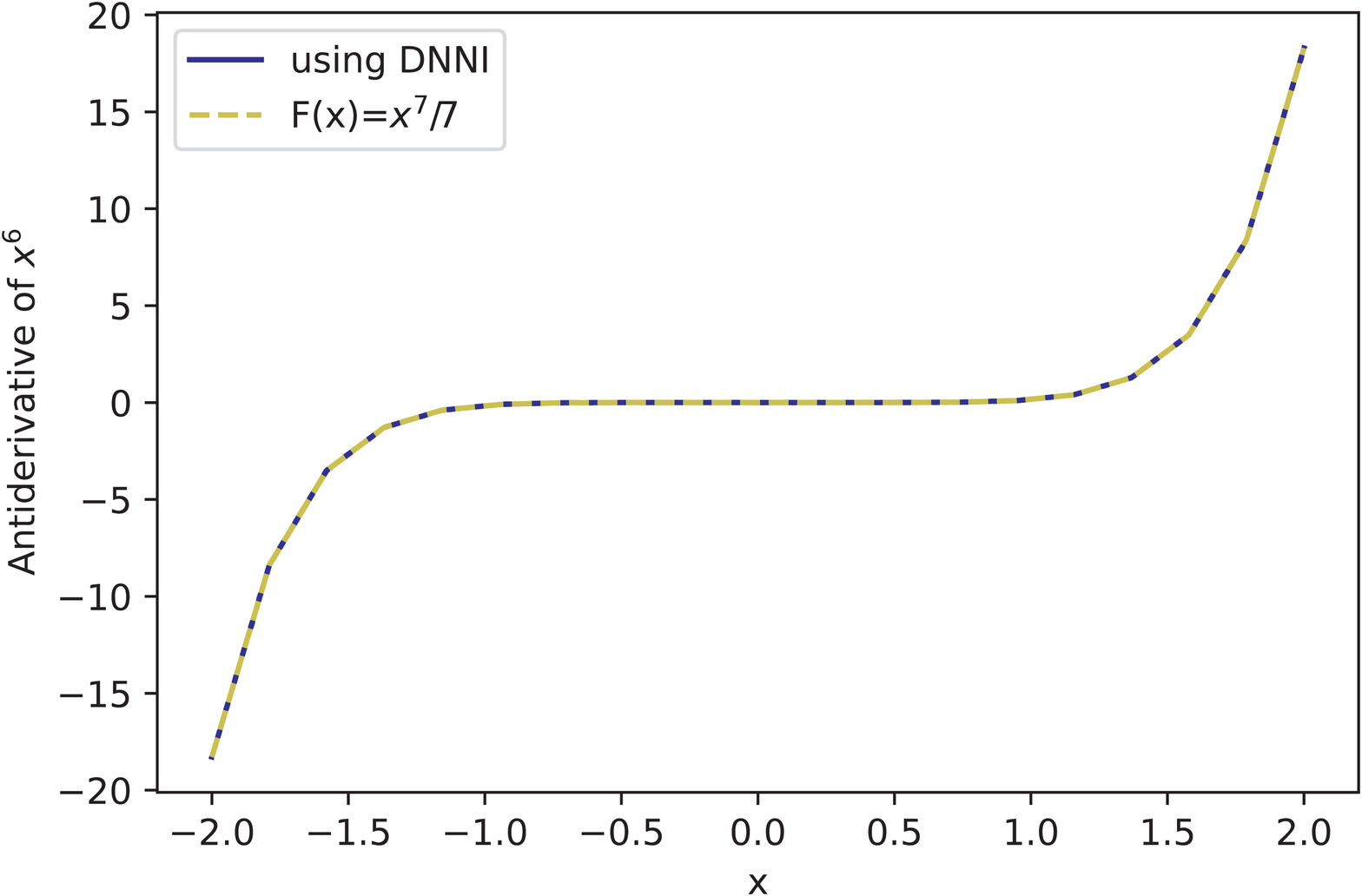}
    \caption{Comparison of the DNNI anti-derivative and theoretical anti-derivative of $x^6$. It can be observed that the DNNI anti-derivative perfectly overlaps with the exact anti-derivative.}
    \label{fig:x^6}
\end{figure}

Case 2: \begin{equation}
    \int \sqrt{1+x^2} \,dx = \frac{\sinh{x}+x(x^2+1)}{2} + c
\end{equation}
\begin{figure}[h]
    \centering
    \includegraphics[width=7.5cm,height=5.6cm]{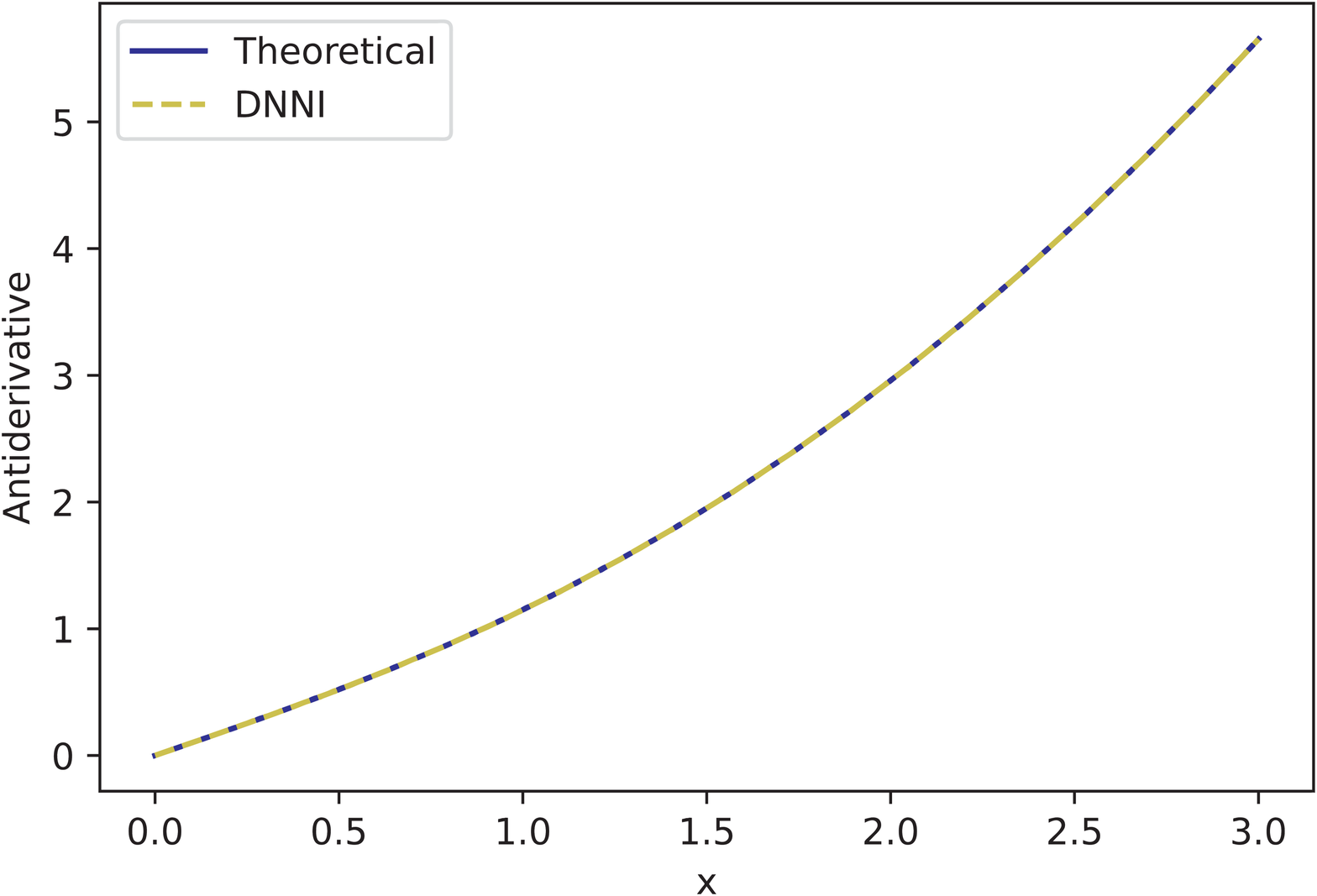}
    \caption{Comparison of the DNNI and theoretical anti-derivative of $\sqrt{1+x^2}$}
    \label{fig:root(1+x^2)}
\end{figure}\\
Case 3: \begin{equation}
    \int \cos{x} \,dx = \sin{x} + c
\end{equation}
\begin{figure}[h]
    \centering
    \includegraphics[width=7.5cm,height=5.6cm]{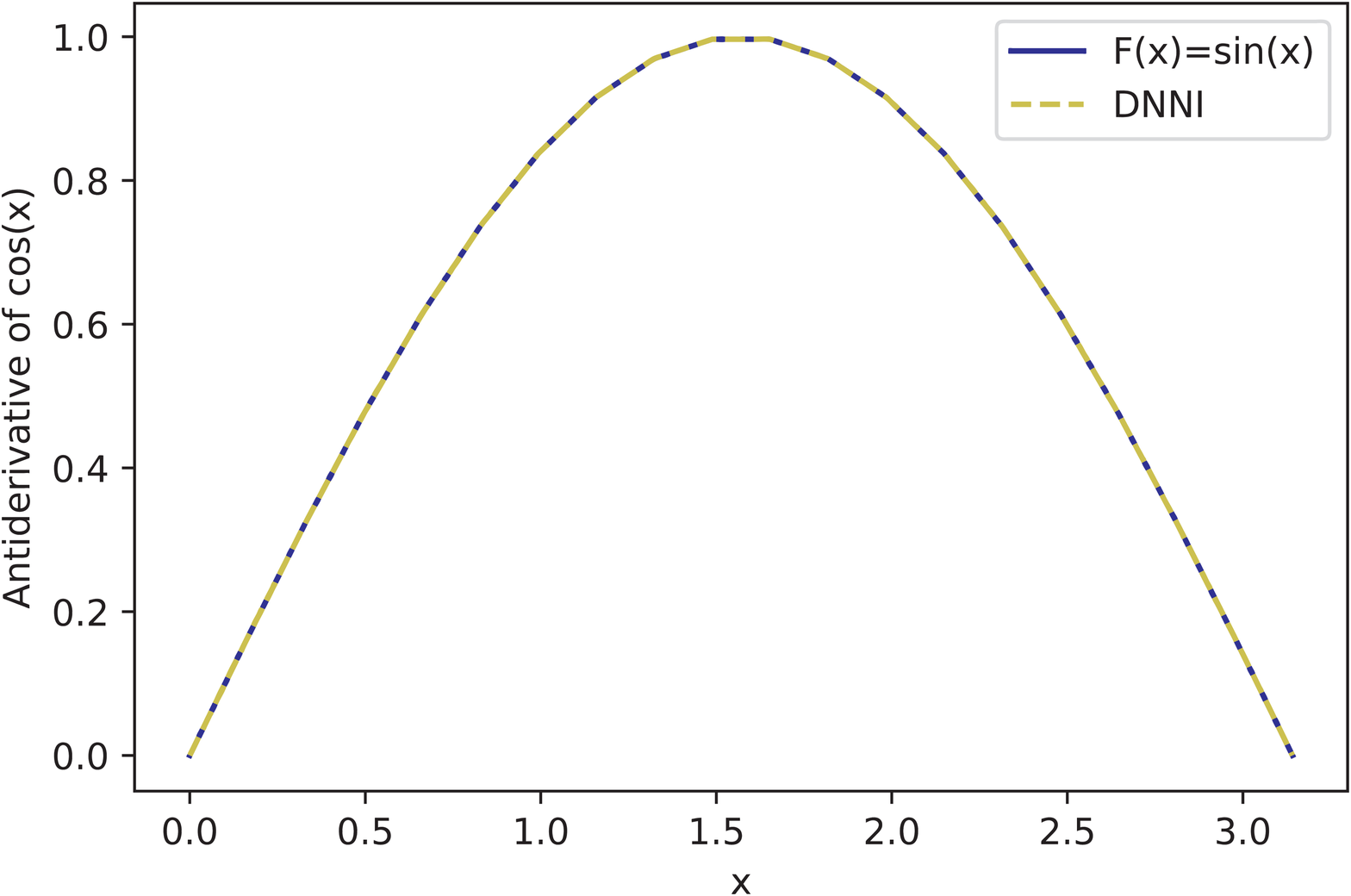}
    \caption{Comparison of the DNNI and theoretical anti-derivative of $cos(x)$}
    \label{fig:cos}
\end{figure}
\newpage
\subsection{Complex Integrals}
Anti-derivatives are very hard to obtain analytically, even if their closed form exists. Applying deep learning in symbolic  mathematics\cite{lample2019deep} has proven helpful in finding primitives of complex integrands. This subsection shows the use of DNNI for such integrals.\\ \\
Case 4:

\begin{equation}\label{com1}
    \int \frac{16x^3-42x^2+2x}{\sqrt{-16x^8+112x^7-204x^6+28x^5-x^4+1}} \,dx  = \sin^{-1}\left( 4x^4-14x^3+x^2 \right) + c
\end{equation}
\begin{figure}[h]
	\centering
	\includegraphics[width=8cm,height=6cm]{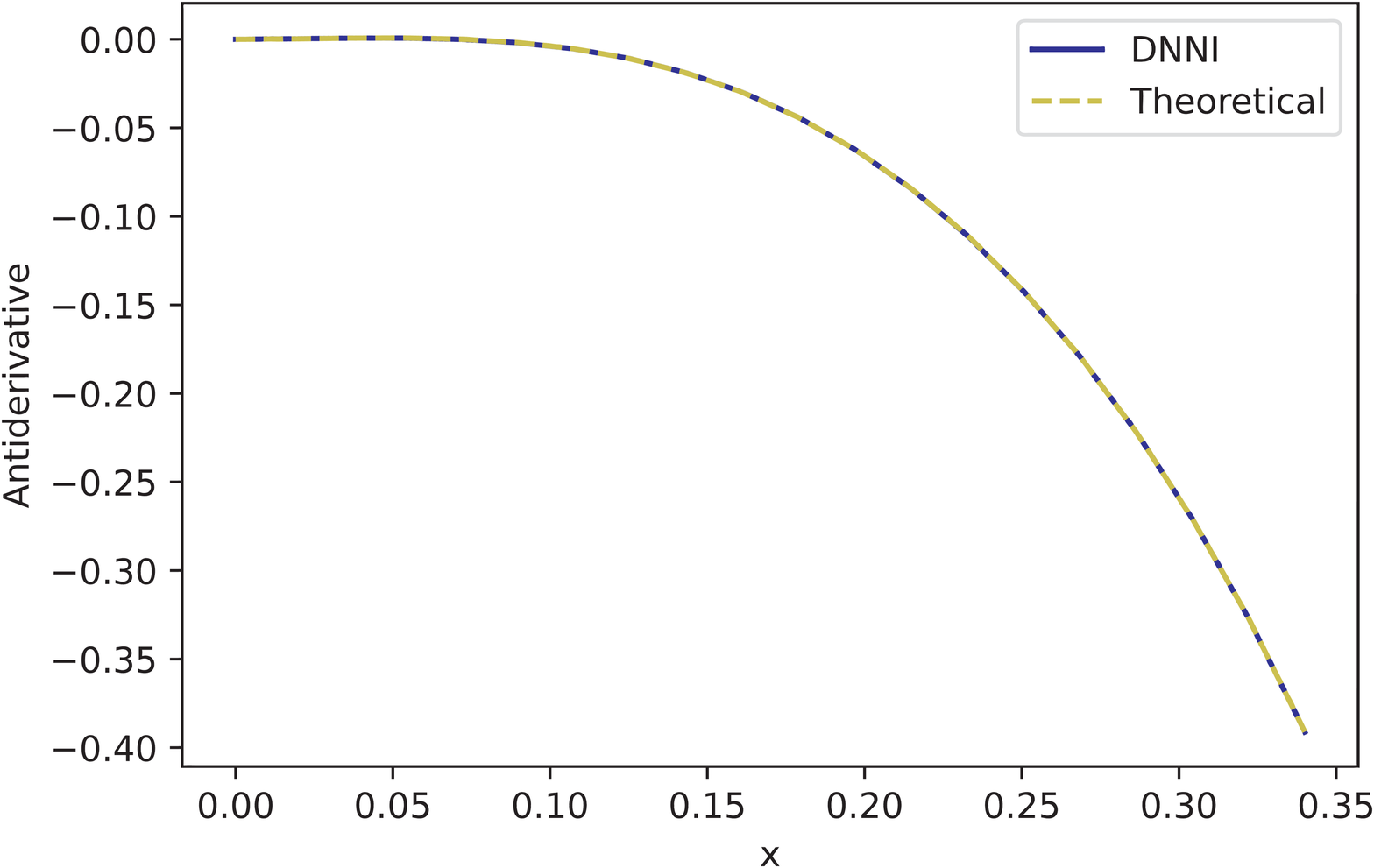}
	\caption{Comparison of the DNNI and theoretical anti-derivative for the integral shown in case 4}
	\label{fig:tough2}
\end{figure}
\\ \\ Case 5\cite{bronstein1998symbolic}:
\begin{equation}\label{com2}
	\int \frac{x^2+2x+1+(3x+1)\sqrt{x+\log(x)}}{x\sqrt{x+\log(x)}(x+\sqrt{x+\log(x)})} \,dx = 2\left(\sqrt{x+\log(x)}+(\log (x+\sqrt{x+\log(x)}) \right) + c
\end{equation}
\begin{figure}[h]
    \centering
    \includegraphics[width=7.5cm,height=5.6cm]{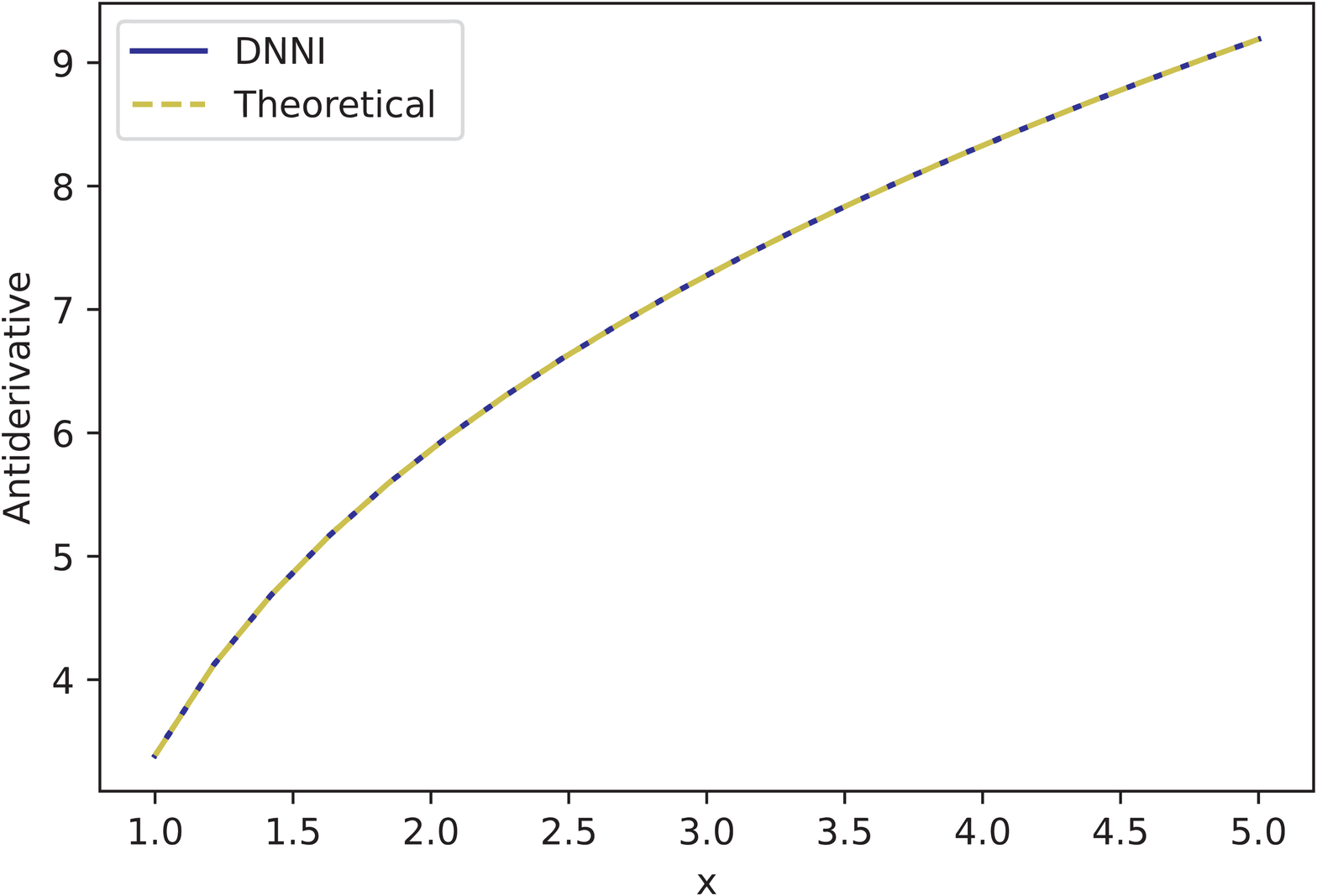}
    \caption{Comparison of the DNNI and theoretical anti-derivative for the integral shown in Case 5}
    \label{fig:tough1}
\end{figure}

\newpage
\subsection{Non-Elementary Integrals}
DNNI is very useful in the case of non-elementary integrals. Though there are several numerical techniques for definite integrals, DNNI can be used to plot the primitive, which is computationally expensive using other numerical techniques due to repeated integrations.
\\ \\ Case 6: 

\begin{figure}[h]
    \centering
    \includegraphics[width=7.5cm,height=5.6cm]{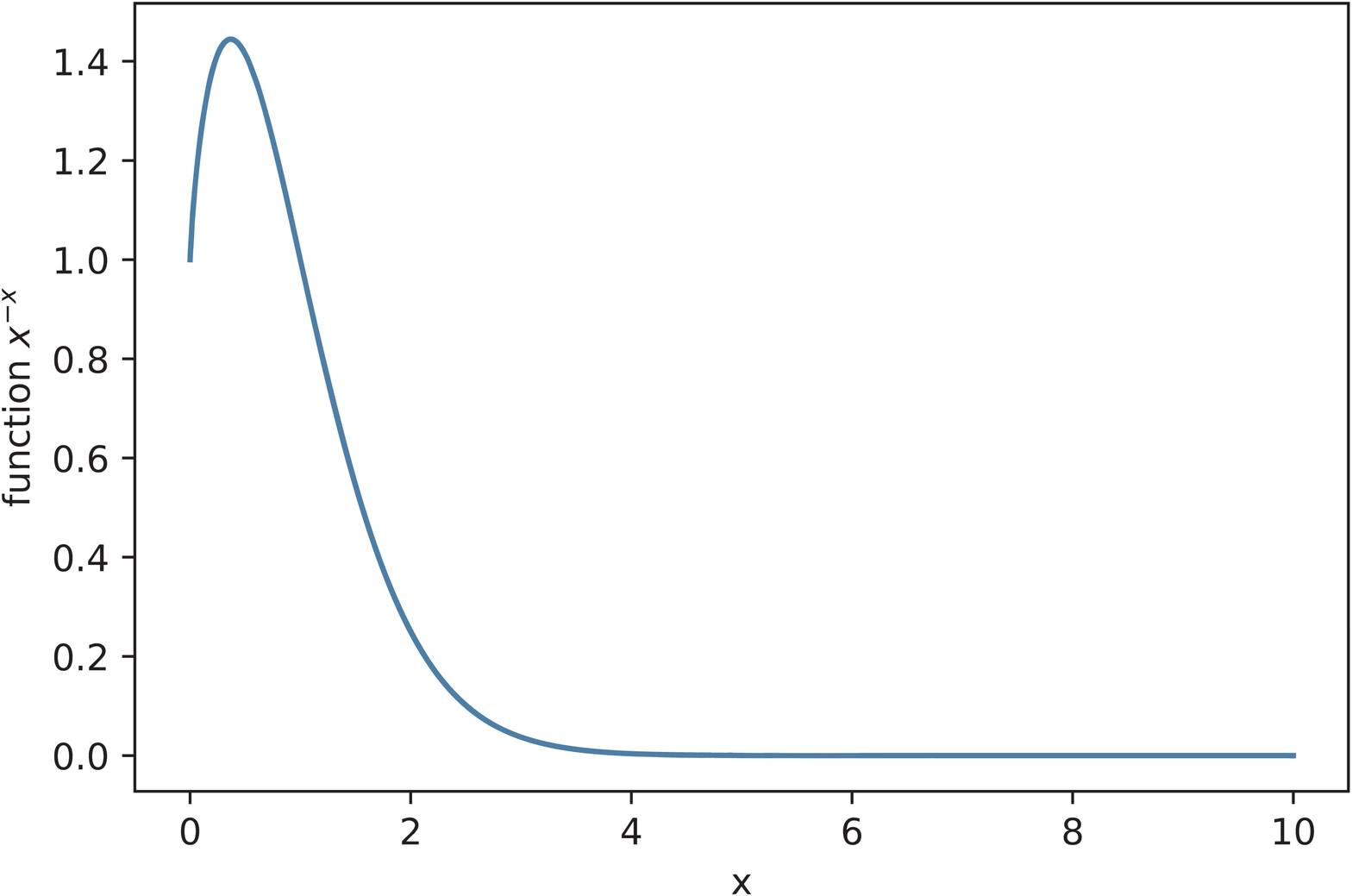}
    \caption{The function $x^{-x}$ which has a maximum value at $1/e$ and then decreases asymptotically}
    \label{fig:x^x}
\end{figure}
The variation of the function $x^{-x}$ is shown in figure \ref{fig:x^x} and its anti-derivative is shown in figure \ref{fig:integral x^x}.
\begin{figure}[h]
    \centering
    \includegraphics[width=7.5cm,height=5.6cm]{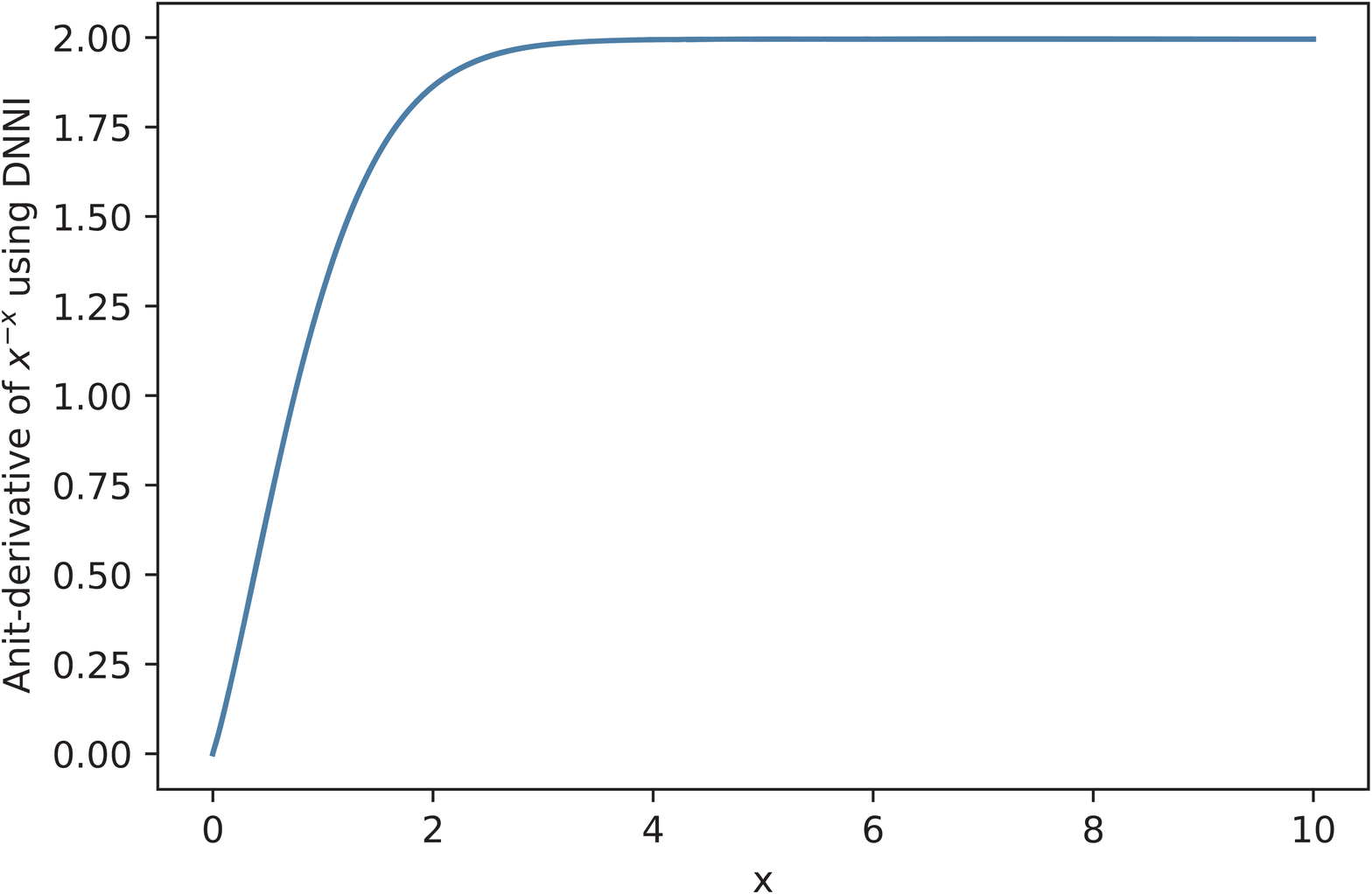}
    \caption{Anti-Derivative of $x^{-x}$ given by $\int_0^x t^{-t} \,dt$ where the lower limit is set to zero for eliminating the constant of integration.}
    \label{fig:integral x^x}
\end{figure}
 \\The identity called a Sophomore's dream is 
\begin{equation}
    \int_0^1 t^{-t} \,dt= \sum_{n=1}^{\infty}n^{-n} = 1.291285997
\end{equation}
which can be verified using DNNI. Also the integral 
\begin{equation}\label{ne2}
    \int_0^\infty t^{-t} \,dt= 1.99545596 
\end{equation}
can be obtained to any desirable accuracy by changing the number of epochs and depth of the neural network. It is obtained to an error of 0.01\% by using 4 hidden layers with 10 nodes each.
\\ \\ Case 6: Elliptic integrals
\\ An elliptic integral is expressed as \begin{equation}
	E(x) = \int_c^x f(t,\sqrt{P(t)})\,dt
\end{equation}
where $f$ is a rational function, and $P$ is a polynomial of degree 3 or 4.
Something as simple as finding the perimeter of an ellipse requires solving a non-elementary integral. The integral to find the perimeter of an ellipse is expressed as
\begin{equation}\label{19}
	Perimeter =\int_0^{\pi/2} 4a\sqrt{1-e^2sin^2(x)}\,dx
\end{equation}
where 'a' and 'b' are major and minor axis lengths and 'e' is the ellipse's eccentricity.\\ \begin{center}
	\begin{tabular}{ c c c c} 
		\hline\hline
		Sl.No. & a & b & Perimeter using Naive DNNI   \\ 
		\hline\hline
		1& 8 & 7 &  47.17621557  \\ 
		\hline
		2 & 2& 1 & 9.68845137 \\ 
		\hline3   & 10 & 5 & 48.44226631 \\
		\hline
		4   & 5 & 1 & 21.01007226 \\\hline\hline
	\end{tabular}
\end{center}
The values obtained have a maximum error of 0.0002\%  using a two hidden layer neural network with ten nodes each. This case is computed using DNNI based on a single variable similar to definite integral calculations. Another approach can be to obtain an approximate closed-form formula based on parameters 'a' and 'b .'This pioneering technique can estimate the closed-form expression of any integral based on several parameters. Further details of this are mentioned in section \ref{param}.

\subsection{Oscillatory Integrals}
Functions of the form $f(x)sin(\frac{\omega}{x^k})$ and $f(x)cos(\frac{\omega}{x^k})$ where $f$ is a continuous and smooth function, $\omega$ and $k$ are real numbers, are highly oscillatory. Integrating such functions is very challenging using common numerical techniques. Mathematicians have developed special techniques like using Haar wavelets and hybrid functions\cite{shivaram2016numerical} to counter such functions. This subsection shows that DNNI can handle even highly oscillatory integrals.  
\\Case 8:
\begin{equation}
    \int_{0}^1 x\ sin(\frac{1}{x^{10}}) \,dx = 0.060665
\end{equation}
\clearpage

\begin{figure}[h]
	\centering
	\begin{minipage}{.5\textwidth}
		\includegraphics[width=\textwidth,height=0.6\textwidth]{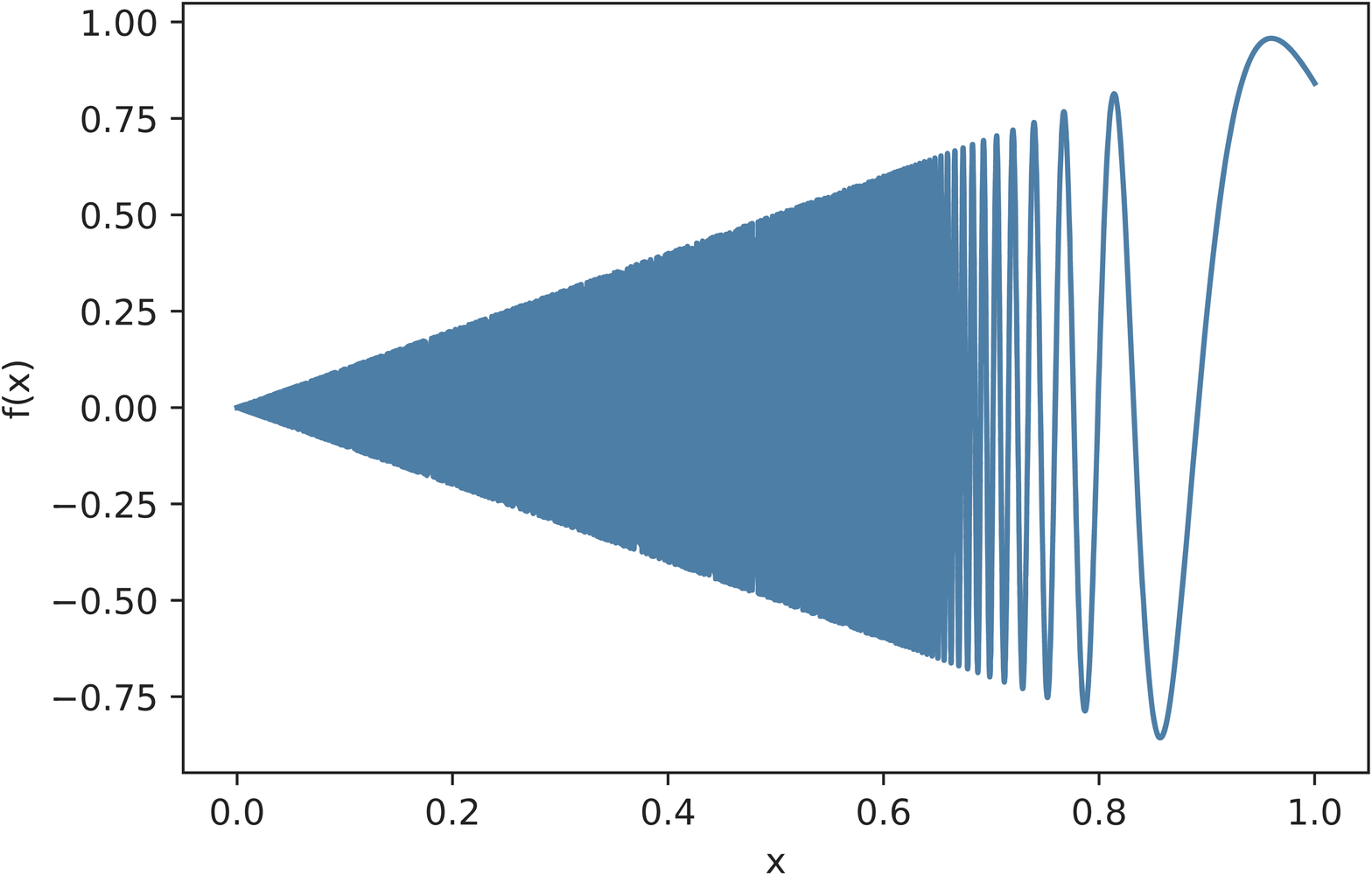}
	\end{minipage}%
	\begin{minipage}{.5\textwidth}
		\includegraphics[width=\textwidth,height=0.6\textwidth]{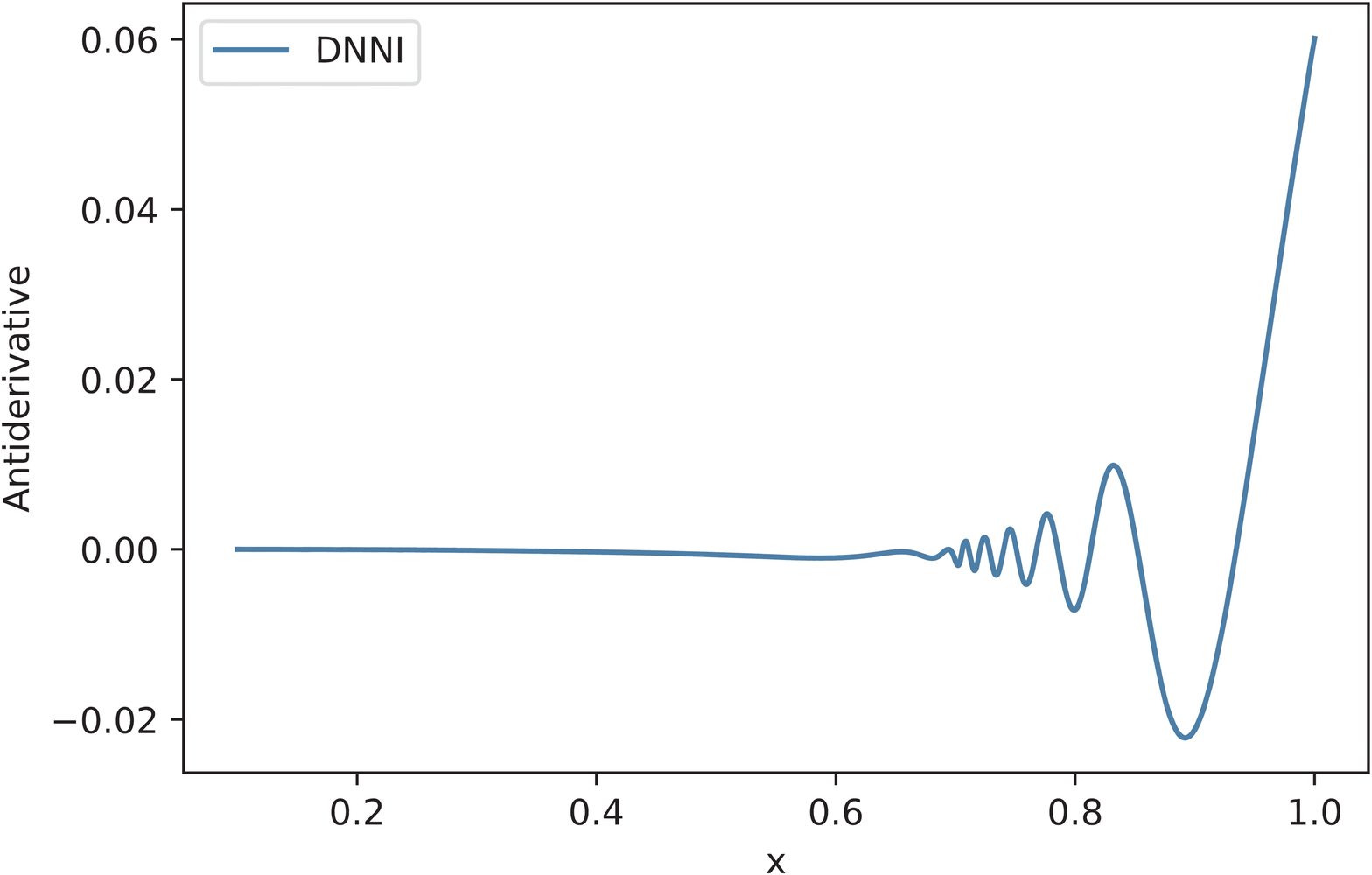}
	\end{minipage}
	\caption{Function: $x\ sin(\frac{1}{x^{10}})$ and its DNNI anti-derivative}
	\label{fig: Osc1}
\end{figure}
Since DNNI approximates the primitive, it is bound to give the correct values of the definite integral on applying the limits. It also gives a closed-form approximation of the anti-derivative, which can be used as required.
\begin{table}[h]
\centering
    \begin{tabular}{c c c}
         \hline\hline
          Method & Value & Error(\%) \\ 
         \hline\hline
         Simpsons 1/3rd(500 points)& 0.0622533209 & 2.61818  \\
         \hline
          Simpsons 3/8th(500 points)& 0.070571762 & 16.33028 \\
         \hline
          Simpsons 1/3rd(1 million points)& 0.0606172467 & 0.07872 \\
          \hline
          Simpsons 3/8th(1 million points)& 0.0605936399 &0.11763 \\ 
         \hline
          Clenshaw-Curtis method(scipy library)& 0.060524  & 0.232424 \\ 
         \hline
          Global Adaptive Quadrature(Matlab default)&0.0605935019 & 0.117857 \\ 
         \hline
          DNNI& 0.06067391 & 0.01469 \\ 
         \hline \hline
    \end{tabular}
	\caption{\label{tab:osc1}A comparison between DNNI and several common numerical techniques.  The computation time of other methods is lower than DNNI, but the number of points used is much higher. DNNI gives the most accurate result among all.}
\end{table}
\\Case 9:
\begin{equation}
    \int_0^1 \frac{1}{x+1}\ sin(\frac{1}{x}) \,dx = 0.28749061
\end{equation}

\begin{figure}[h]
	\centering
	\begin{minipage}{.5\textwidth}
		\includegraphics[width=\textwidth,height=0.6\textwidth]{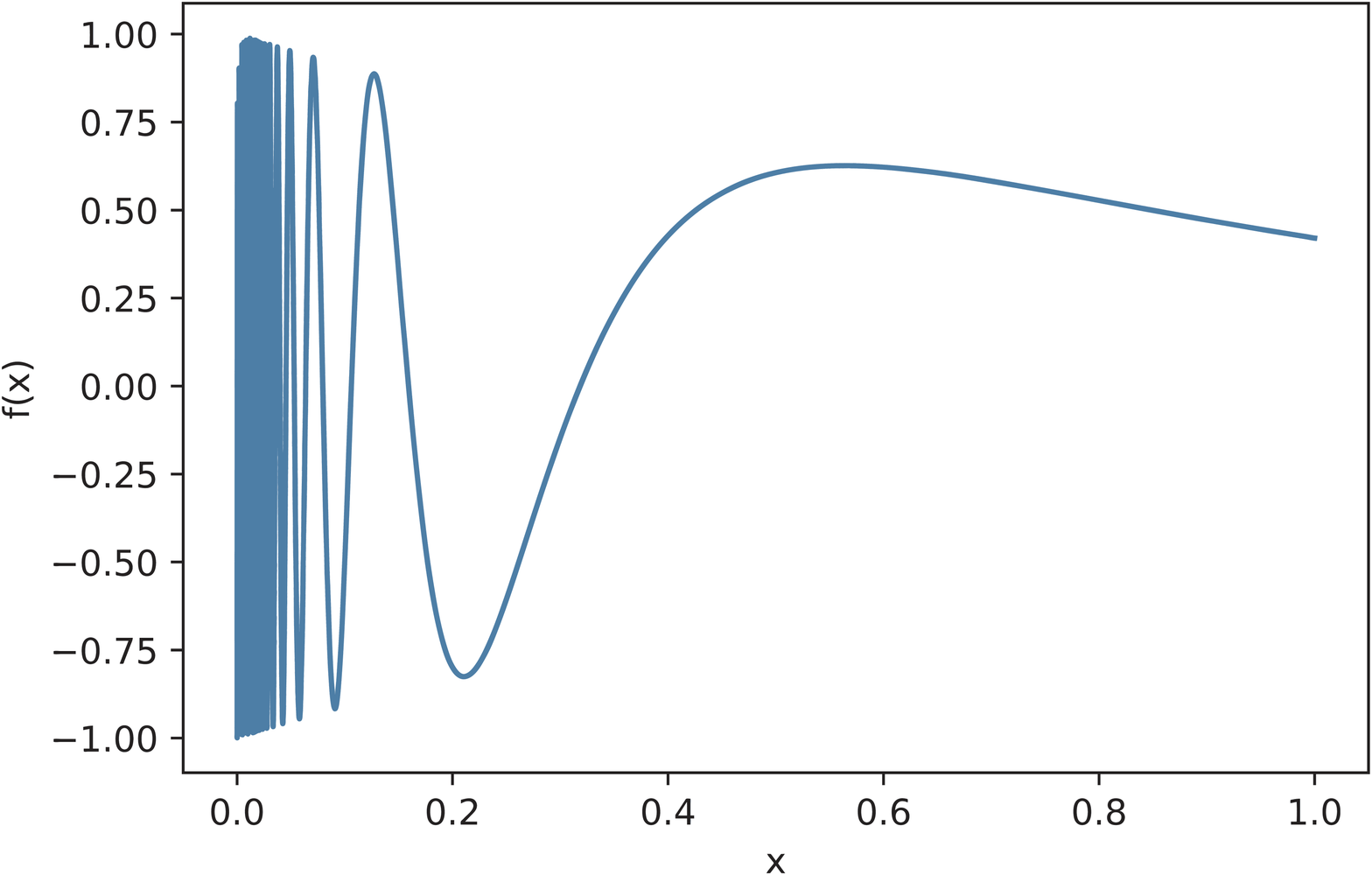}
	\end{minipage}%
	\begin{minipage}{.5\textwidth}
		\includegraphics[width=\textwidth,height=0.6\textwidth]{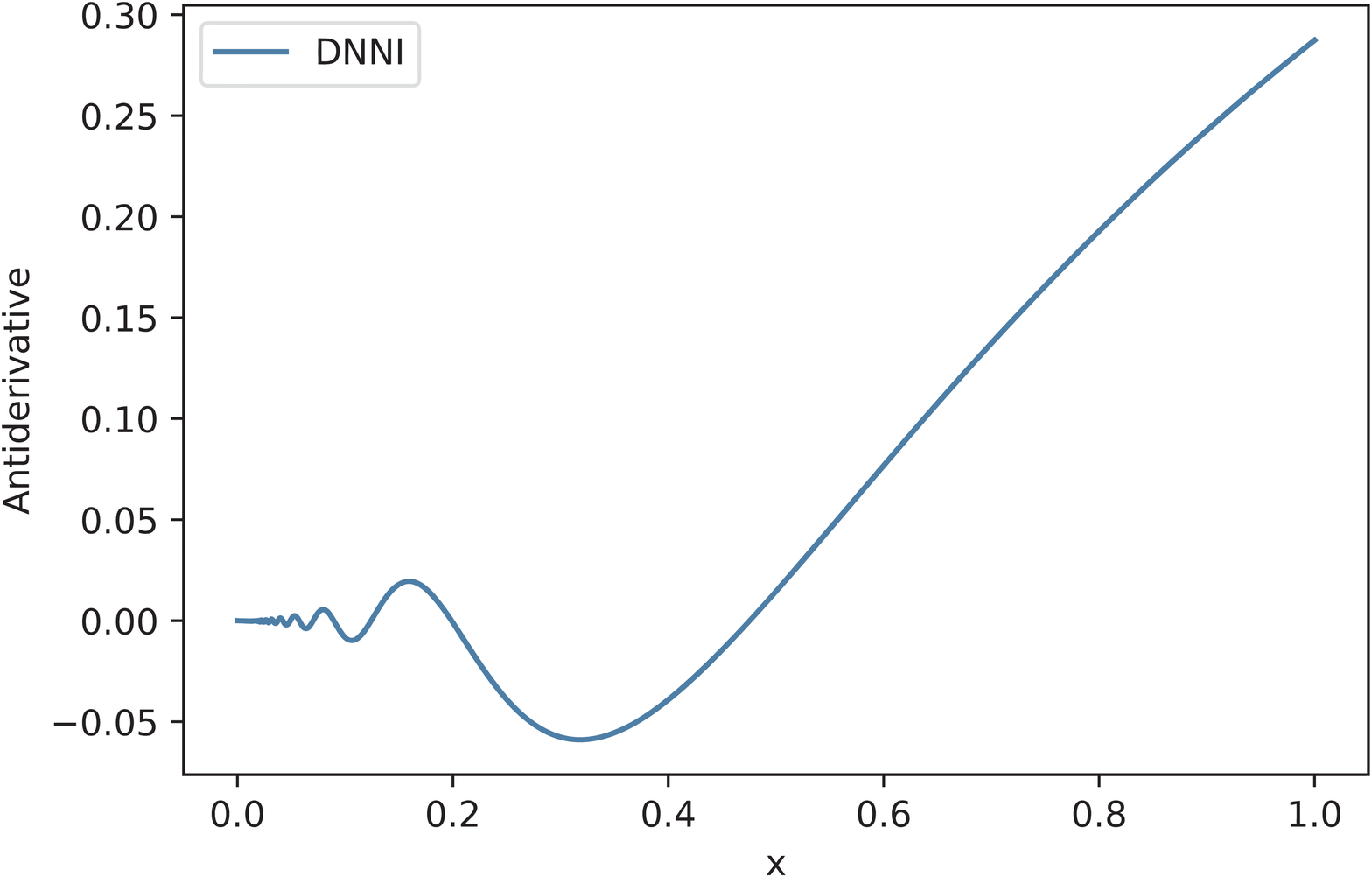}
	\end{minipage}
	\caption{ The function: $\frac{1}{x+1}\ sin(\frac{1}{x})$ and its DNNI anti-derivative}
	\label{fig: Osc2}
\end{figure}
\begin{table}[h]
\centering
    \begin{tabular}{c c c}
         \hline \hline
          Method & Value & Error(\%) \\ [0.5ex] 
         \hline\hline
         Simpsons 1/3rd(500 points)& 0.28603691 & 0.50565  \\
         \hline
          Simpsons 3/8th(500 points)& 0.28345895 & 1.402362 \\
         \hline
          Simpsons 1/3rd(1 million points)& 0.28751143 & 0.007242 \\
          \hline
          Simpsons 3/8th(1 million points)& 0.28750075 & 0.003527 \\ 
          
         \hline
          Clenshaw-Curtis method(scipy library)& 0.285857 & 0.568230 \\ 
         \hline
         
         Global Adaptive Quadrature(Matlab default)&0.28749060 & Exact \\
         \hline
          DNNI& 0.28730544 & 0.064409 \\ 
         \hline \hline
    \end{tabular}
	\caption{\label{tab:osc2} A comparison between DNNI and several common numerical techniques. Since this function is less oscillatory than case 8, the numerical techniques perform relatively better.}
\end{table}

\newpage
\subsection{Error Analysis}
The accuracy of DNNI increases with the number of points taken for training the neural network. The l2 norm from the theoretical solution decreases asymptotically with the number of training points showing some local fluctuations, which decay out on average. In this paper, the learning rate is of the order $10^{-2}$ and is reduced every one-fifth of the training steps.  Increasing the number of epochs with decreasing learning rates also decreases the l2 norm asymptotically. We can tune all these parameters accordingly to obtain the best model.  
\begin{table}[h]
\centering
\begin{tabular}{c c c}
	\hline \hline
	Parameter & Simple Integral & Complex Integral \\ [0.5ex] 
	\hline\hline
	No. of training points & 20-100 & 1000-5000  \\
	\hline
	Depth of the Neural Network& 2-4 layers & 4-8 layers \\
	\hline
	No. of nodes in each layer& 5-10 & 10-20 \\
	\hline
	No. of epoches& 1000-10000 & 10000-50000 \\ 
	\hline \hline
\end{tabular}
\caption{We suggest using the above parameters to train the DNNI model based on the complexity of the problem. A meshgrid has to be generated for the integrals based on several parameters to use all combinations in the training data.}
\end{table}
\begin{figure}[h]
	\centering
	\begin{minipage}{.5\textwidth}
		\centering
		\includegraphics[width=\linewidth,height=0.6\linewidth]{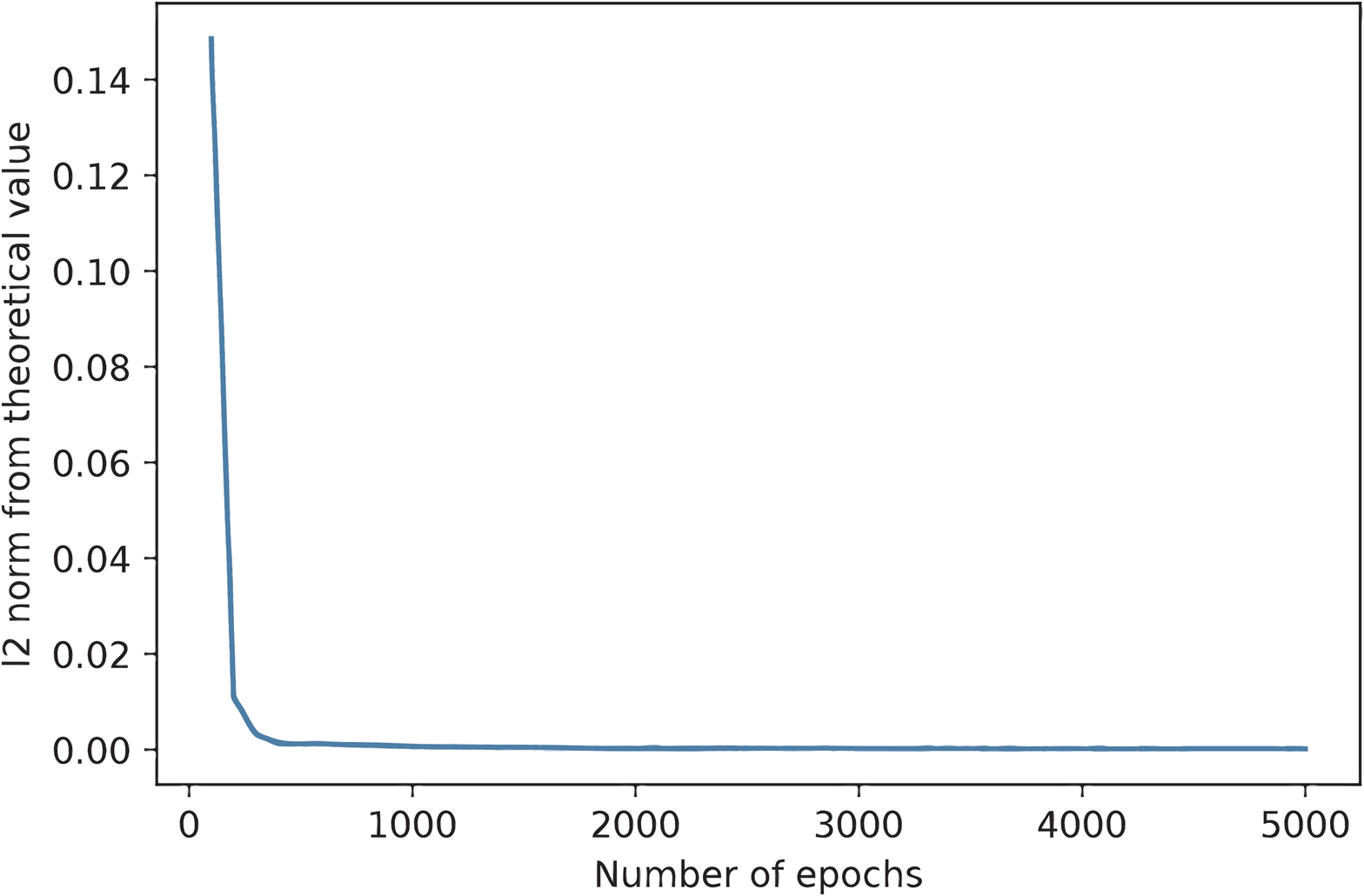}
		\label{fig:e1}
	\end{minipage}%
	\begin{minipage}{.5\textwidth}
		\centering
		\includegraphics[width=\linewidth,height=0.6\linewidth]{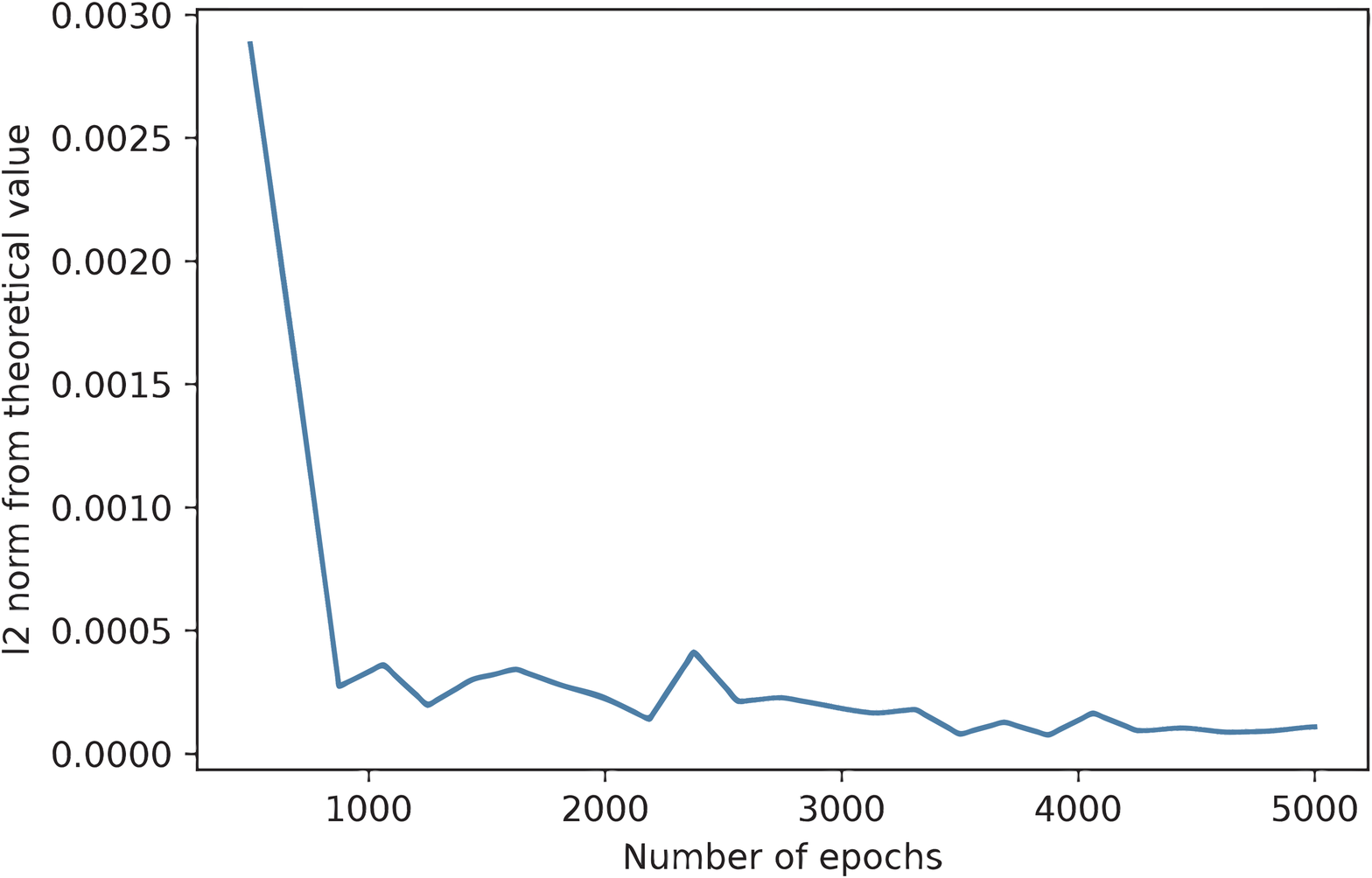}
		\label{fig:e2}
	\end{minipage}
	\caption{The variation of L2 norm with respect to theoretical anti-derivative for the integrals in equation \ref{com1} and \ref{ne2} respectively with increasing epochs.}
	\label{fig:epoch}
\end{figure} 

\begin{figure}[h]
	\centering
	\begin{minipage}{.5\textwidth}
		\centering
		\includegraphics[width=\linewidth]{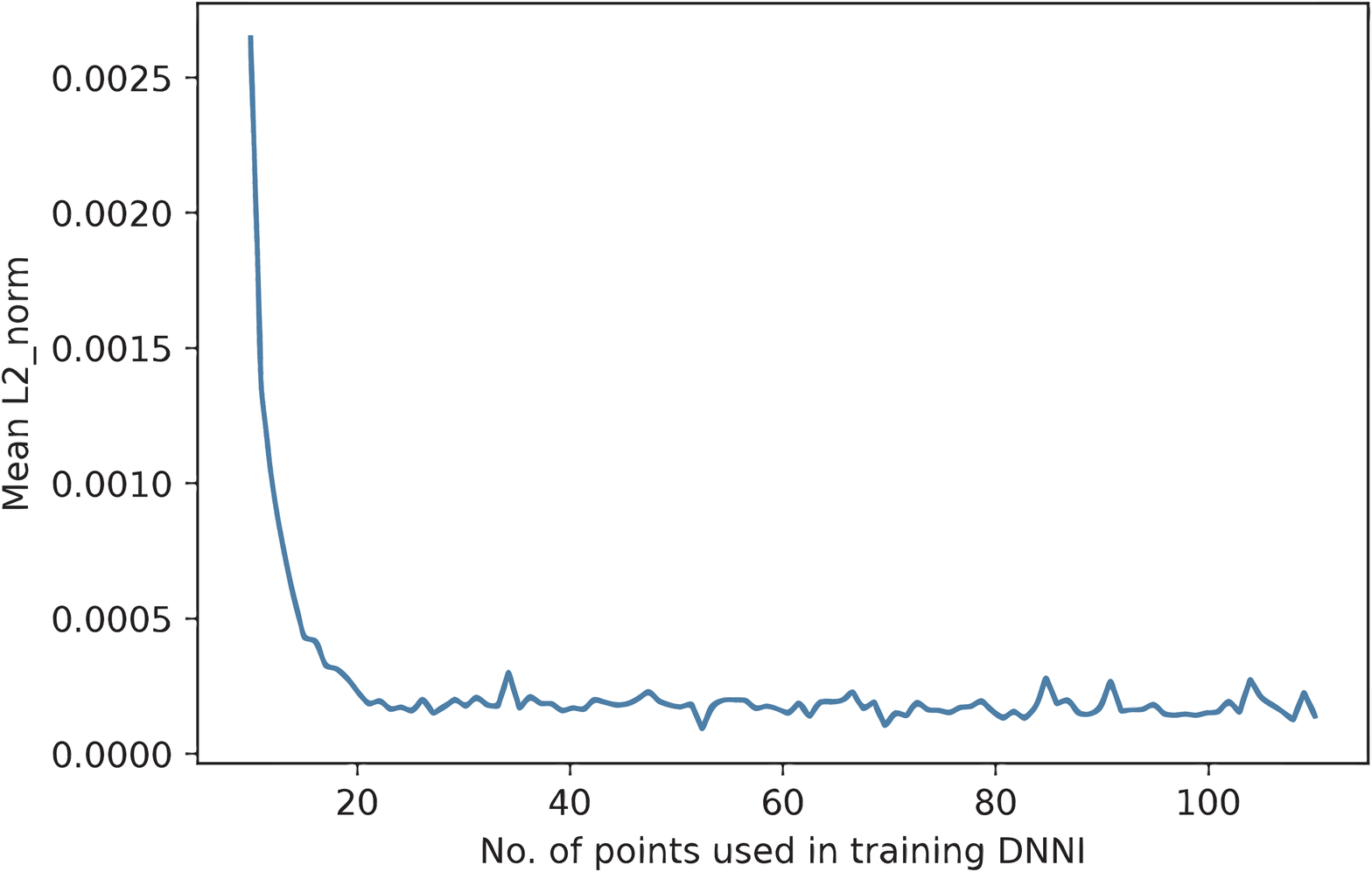}
		\label{fig:sub1}
	\end{minipage}%
	\begin{minipage}{.5\textwidth}
		\centering
		\includegraphics[width=\linewidth]{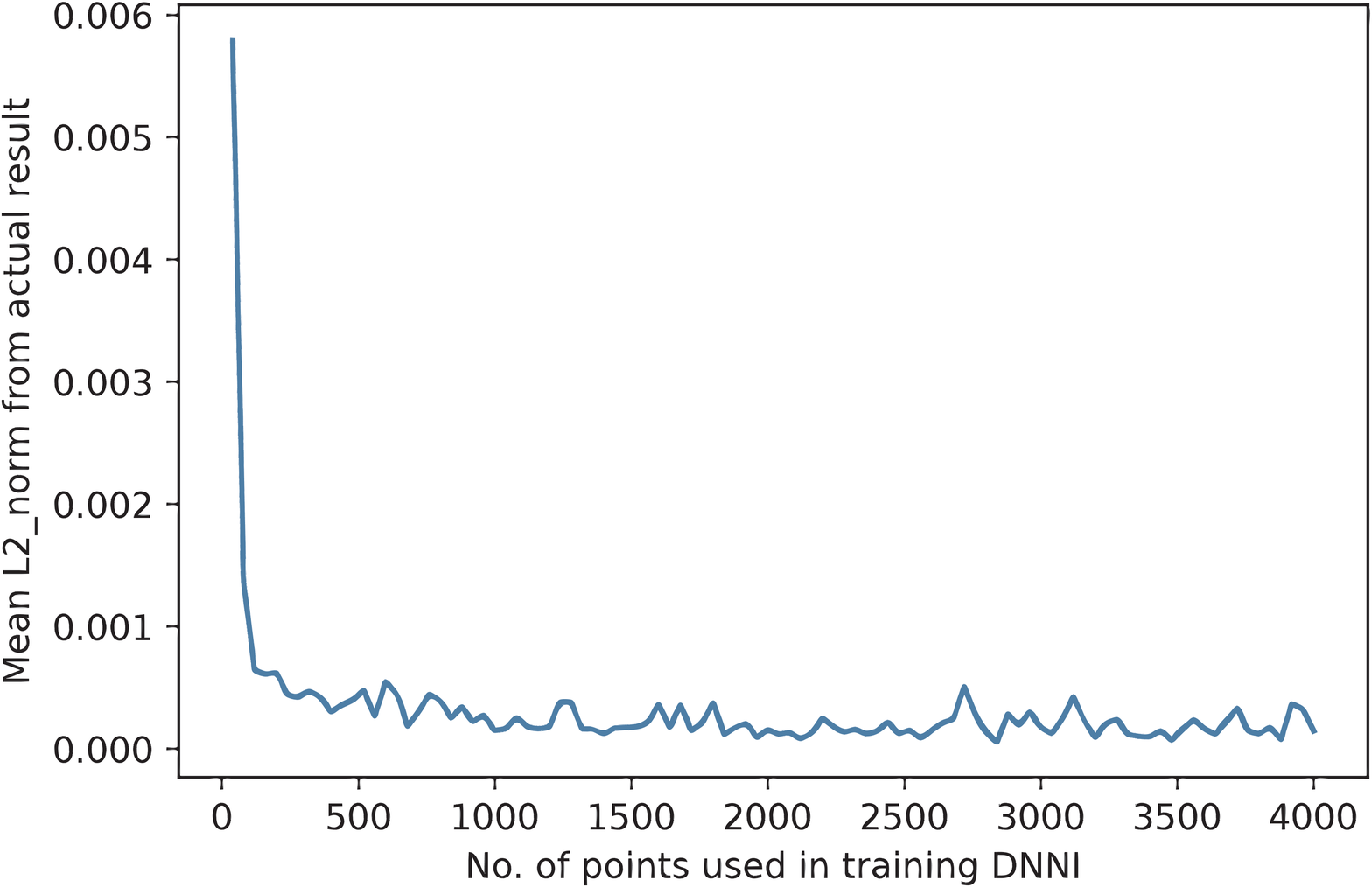}
		\label{fig:sub2}
	\end{minipage}
	\caption{The variation of L2 norm with respect to theoretical anti-derivative for the integrals in equation \ref{com1} and \ref{ne2} respectively on increasing the number of training points.}
	\label{fig:error1}
\end{figure}  

%\clearpage
%\newpage

\newpage
\section{ Applications}
The applications of DNNI that we have come up with are based on either the requirement of anti-derivatives of non-elementary integrals or increasing the computational speed of the algorithms having repeated integrals.
\subsection{Parametric integrals}\label{param}
This subsection shows the use of DNNI in obtaining anti-derivatives of integrands having several parameters. The closed form solution is required in such cases to use the integral further and study the effects of various parameters. Mathematicians have successfully developed many closed-form approximations for popular non-elementary integrals like the perimeter of an ellipse and the Fermi-Dirac integral. Using these test cases, we claim that DNNI can be used to obtain the closed form approximation of any given integral. Though there is no substitute for theoretical analysis, DNNI can effectively give a relatively quick solution in very complex cases. Also, DNNI is a single method that works for all types of integrands.\\
Case 10:\\
The equation \ref{19} shows the integral for the perimeter of
an ellipse. This integral depends on the semi-major axis 'a' and the semi-minor axis 'b.' 
\begin{equation}
	P(a,b) = 4 \int_0^{\pi/2} \sqrt{a^2-(a^2-b^2)sin^2(\theta)} \,d\theta,\ \ \ \ a,b\in\mathbb{R}
\end{equation}
For any given 'a' and 'b,' the value can be found easily using numerical definite integral techniques, but a closed form expression of P(a,b) is often required. The above integral is very common and has been approximated by several mathematicians, including the famous Ramanujan's formula:
\begin{equation}
	P(a,b) \approx a\left(3(a+b)- \sqrt{\frac{3a+b}{a+3b}}\right) 
\end{equation}
DNNI can be used to obtain an approximate closed-form anti-derivate of the expression. The inputs to the neural network have to be a flattened meshgrid of 3 parameters: $\theta$, $a$, and $b$ in the required domain of interest.
The perimeter can be calculated as:
\begin{equation}
	P(a,b) \approx N(\frac{\pi}{2},a,b)-N(0,a,b) 
\end{equation}
\begin{table}[h]
	\begin{center}
		\begin{tabular}{ c c c c} 
			\hline \hline
			a & b &   $N(\frac{\pi}{2},a,b)-N(0,a,b)$  & Relative Error \\ 
			\hline \hline
			5 & 1 &  21.03439167 & 0.001159 \\ 
			\hline
			6& 1.8 & 26.29762002 & 0.000858 \\ 
			\hline 7 & 2.6 & 31.75970172 &0.000171 \\
			\hline
			8 & 3.4 &  37.28223005 &0.000139 \\\hline
			9 & 4.2&  42.84975621 & 0.000043\\\hline
			10&5& 48.40454929 & 0.0000078\\ \hline \hline
		\end{tabular}
	\end{center}
\caption{DNNI is used to obtain a formula for the perimeter of an ellipse. These errors are slightly higher than the ones obtained in Case 6, which was based on a single variable and fixed parameters. The errors can be further reduced by using deeper neural networks and increasing the number of epochs.}
\end{table}
\newpage
Case 11:\\
The non-relativistic Fermi-Dirac integral is defined as
\begin{equation}\label{fermi}
	F_q(\eta) = \int_0^{\infty} \frac{x^q}{e^{x-\eta}+1} \,dx ,\ \ \ q\geq0,\ \ \ n\in\mathbb{R}
\end{equation}
and the relativistic Fermi-Dirac integral is
\begin{equation}
	F_q(\eta,\beta)=\int_0^{\infty} \frac{x^q\sqrt{1+\beta x/2}}{e^{x-\eta}+1}\,dx,\ \ \ \beta\geq0,\ \ \ q\geq0,\ \ \ n\in\mathbb{R}
\end{equation}
The Fermi-Dirac integral has many applications in nuclear astrophysics and finding the concentration of electrons and holes in a semiconductor. There has been substantial research \cite{gil2022complete,sagar1991gaussian,temme1990uniform,bhagat2003evaluation,mohankumar2016very} on just trying to obtain theoretical and numerical approximates of the relativistic and non-relativistic Fermi-Dirac integral. In the following tables, DNNI is used to get an approximate closed-form expression for the above integrals. The values are compared to those obtained using numerical definite integral techniques. DNNI gives a function that outputs the integral on inputting the parameters. The neural network approximate of the Fermi-Dirac integral can be further integrated, differentiated, and plotted based on the requirements.
\begin{table}[h]
	\begin{center}
		\begin{tabular}{ c c c c} 
			\hline \hline
			q & $\eta$ &   $N(\zeta,q,\eta)-N(0,q,\eta)$  & Relative Error \\ 
			\hline \hline
			0 & -2 &  0.12468052 & 0.017706 \\ 
			\hline
			0.5& -1 & 0.28986771 & 0.002179 \\ 
			\hline 1 & 0 &0.8233424 &0.001064 \\
			\hline
			1.5 & 1 & 2.66133345 &0.000130 \\\hline
			2 & 2& 9.51024877 & 0.000254\\\hline \hline
		\end{tabular}
	\end{center}
	\caption{DNNI is used to obtain a formula for the non-relativistic Fermi-Dirac integral. Here, $\zeta$ is an arbitrarily large number in the given domain. This table shows the relative errors of a few integrals by putting the required parameters in the DNNI function with their definite integral counterpart. The errors can be further reduced by using deeper neural networks and increasing the number of epochs.}
\end{table}
\\The multi-parameter DNNI is computationally expensive compared to a single variable DNNI, but it has extensive use in all fields of science and engineering.
\begin{table}[h]
	\begin{center}
		\begin{tabular}{ c c c c c} 
			\hline \hline
			q & $\eta$ & $\beta$ &  $N(\zeta,q,\eta,\beta)-N(0,q,\eta,\beta)$  & Relative Error \\ 
			\hline \hline
			1 & -1 &  0.5 & 0.41499549 & 0.0000397 \\ 
			\hline
			1.5& 0 & 1 & 1.74834439 & 0.005471 \\ 
			\hline 2 & 1 & 1.5 &7.94319678 &0.001817 \\
			\hline
			2.5 & 2 &2& 38.88427763 &0.004774 \\\hline \hline
		\end{tabular}
	\end{center}
	\caption{DNNI is used to obtain a formula for the relativistic Fermi-Dirac integral. This table shows the relative errors of a few integrals by putting the required parameters in the DNNI anti-derivative with their definite integral counterpart.}
\end{table}

\newpage

\subsection{Cumulative Distribution function}
\begin{figure}[h]
	\centering
	\begin{minipage}{.5\textwidth}
		\includegraphics[width=\textwidth,height=0.8\textwidth]{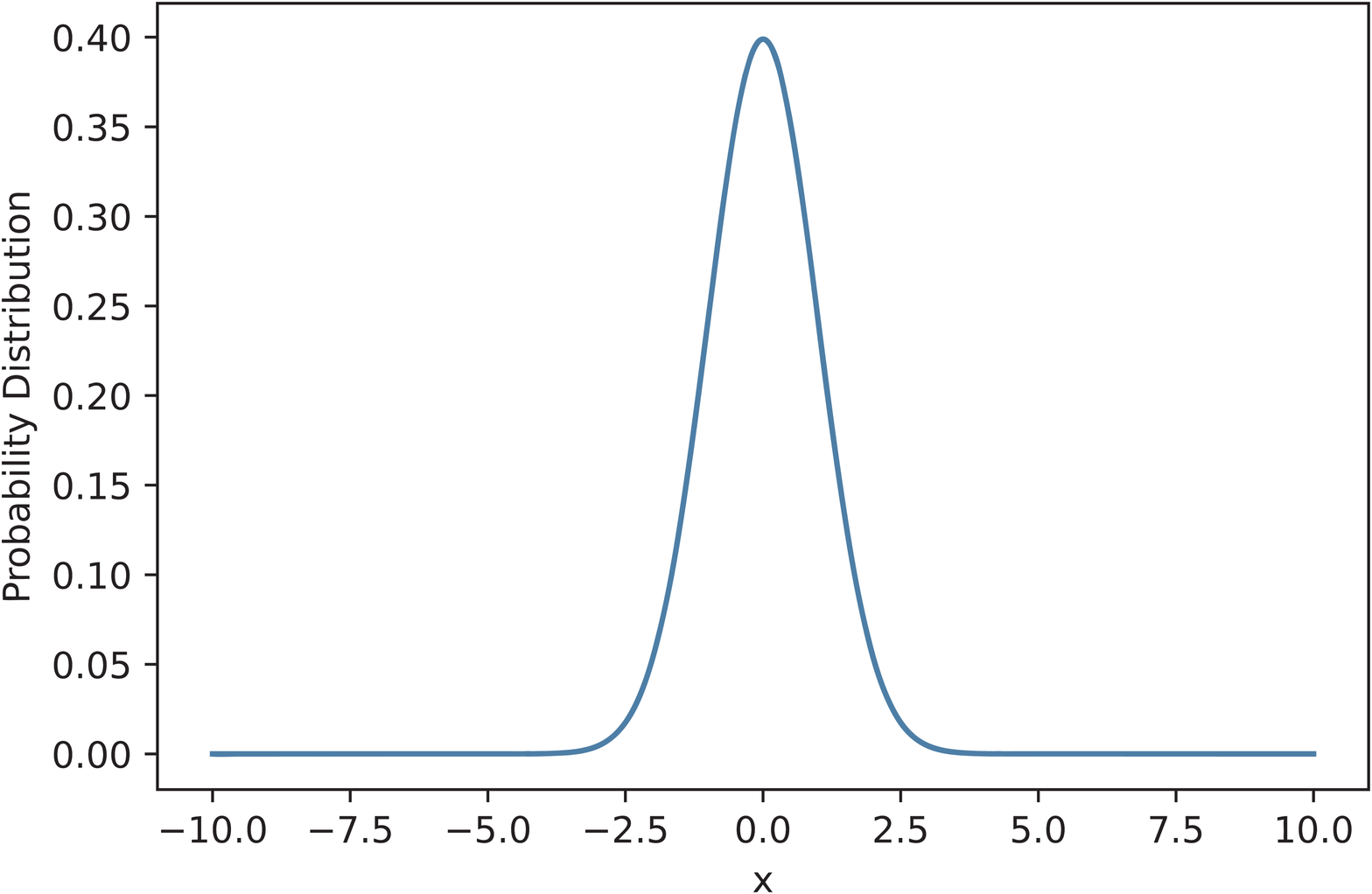}
	\end{minipage}%
	\begin{minipage}{.5\textwidth}
		\includegraphics[width=\textwidth,height=0.8\textwidth]{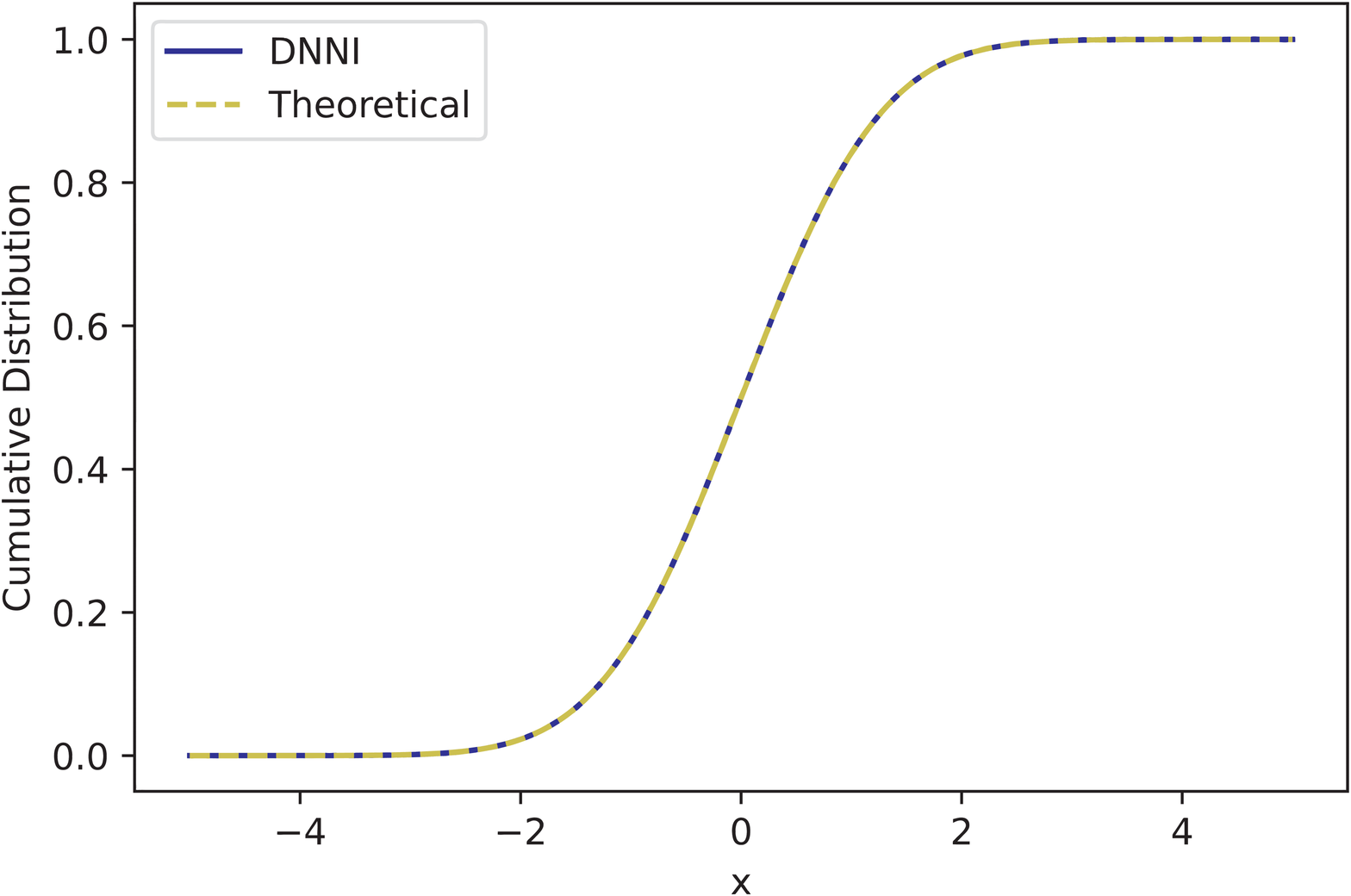}
	\end{minipage}
	\caption{Standard Normal Distribution and its Cumulative Distribution function}
	\label{fig:normal}
\end{figure}

A probability distribution function gives the distribution of the probability of occurrence of an event. For a continuous random variable, it is often represented by a function $f$. The cumulative distribution function of a random variable is given by 
\begin{equation}
	F(X) = \int_{-\infty}^X f(x) \,dx
\end{equation}
The need for an anti-derivative is very critical in this case. DNNI can be a handy method to approximate any cumulative distribution function. It can replace the need for distribution tables or repeated numerical integration. This claim is shown in the following cases: 
\\ \\Case 12:\\
For the probability distribution function 
\begin{equation}
	f(x) = \frac{1}{2\pi}e^{-x^2/2},
\end{equation}
as shown in figure \ref{fig:normal}, the cumulative distribution function is
\begin{equation}
	\frac{1}{2\pi}\int_{-\infty}^x e^{-t^2/2} \,dt = \frac{1}{2}\ (1+erf\left( \frac{x}{\sqrt{2}}\right)) 
\end{equation}

\begin{figure}[h]
	\centering
	\begin{minipage}{.5\textwidth}
		\includegraphics[width=\textwidth,height=0.8\textwidth]{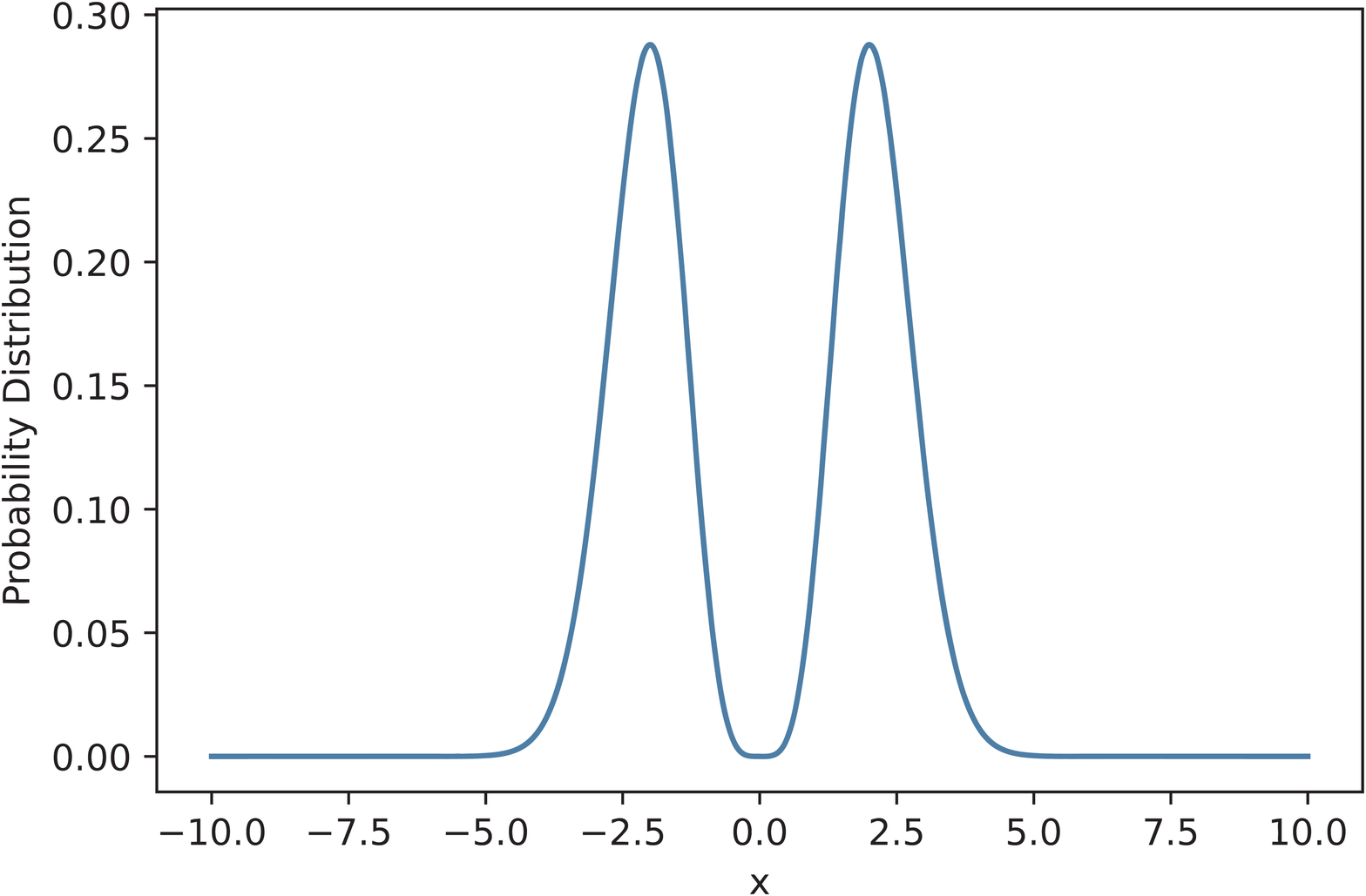}
	\end{minipage}%
	\begin{minipage}{.5\textwidth}
		\includegraphics[width=\textwidth,height=0.8\textwidth]{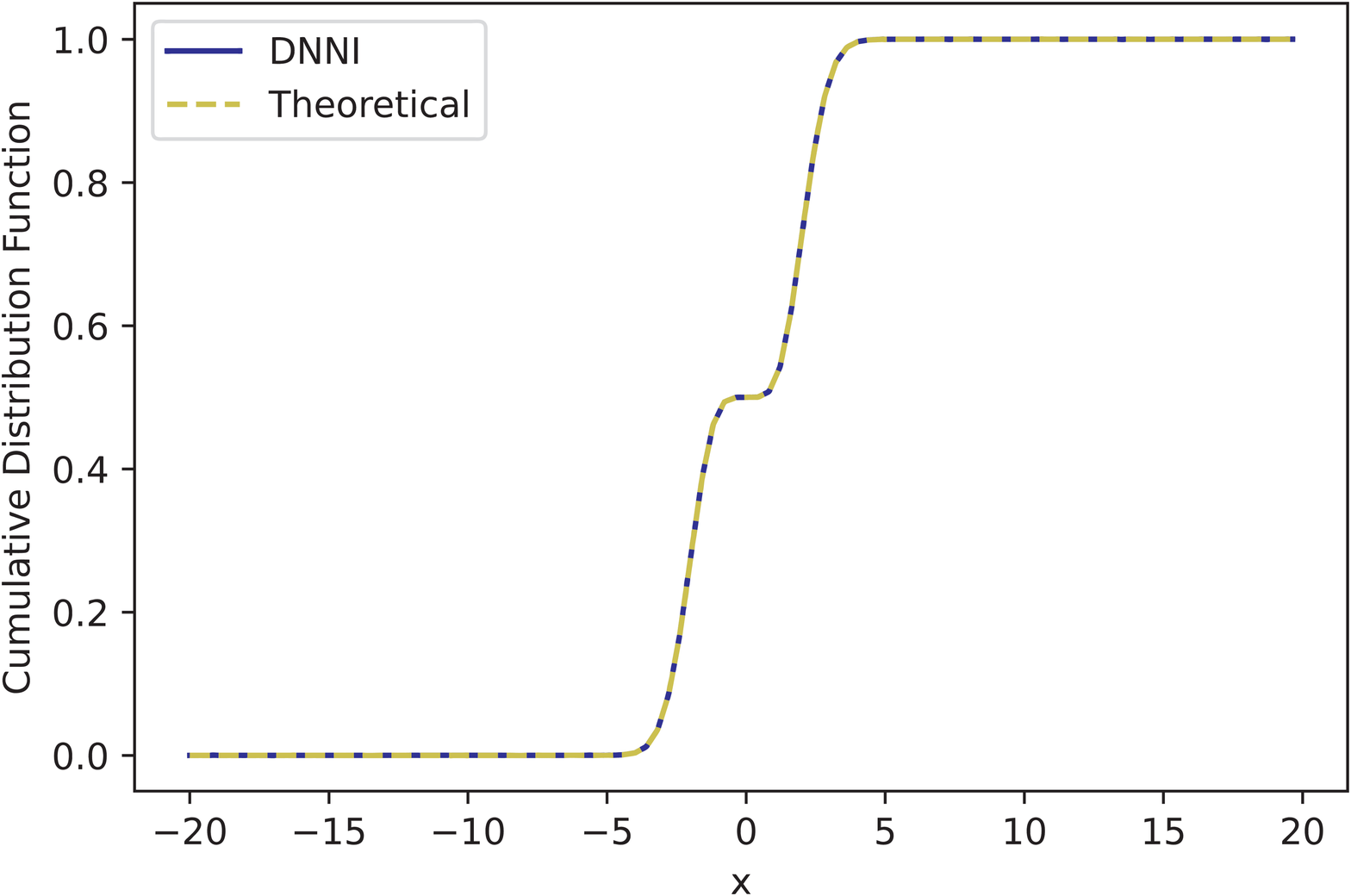}
	\end{minipage}
	\caption{A Bimodal Distribution and its Cumulative Distribution Function}
	\label{fig: A special bimodal Distribution}
\end{figure}

Case 13:
\\
For the probability distribution function  
\begin{equation}
	f(x) = \frac{1}{3\sqrt{2\pi}}x^4e^{-x^2/2},
\end{equation}
as shown in figure \ref{fig: A special bimodal Distribution}, the cumulative distribution function is
\begin{equation}
	\frac{1}{3\sqrt{2\pi}}\int_{-\infty}^x t^4e^{-t^2/2} \,dt 
	= \frac{1}{2}\ (1+erf\left( \frac{x}{\sqrt{2}}\right)) - \frac{1}{3\sqrt{2\pi}}x(x^2+3)e^{-x^2/2}
\end{equation}
\\  Using DNNI, the cumulative distribution function of very complicated probability distributions can also be found. Numerically finding the definite integral values for different limits can be computationally expensive. DNNI can be an up-and-coming solution to this.
\newpage
\subsection{Galerkin Method}
Galerkin Method is a numerical technique to solve differential equations by converting them to weak integral form. The domain is discretized, and the conservation variable is approximated with a basis function. The differential equation is multiplied with a weight function and integrated such that the residual is zero. The basis function's coefficients are calculated by solving sets of linear equations. Galerkin Method is a very general and broad family of methods\cite{giraldo2020introduction} to solve differential equations. We have picked one such method closer to the Finite Element Method.
\\For a differential equation,
\begin{equation}
	\alpha\frac{d^2u(x)}{dx^2}+\beta \frac{du(x)}{dx}+\gamma u(x) = f(x)
\end{equation}
with Dirichlet boundary condition $u(x_0)=u_0$ and Neumann boundary condition $\frac{du}{dx}$=c at $x=x_n$ in the domain $x_0<x<x_n$.
Let $u(x)$ be approximated as
\begin{equation}\label{eqn28}
	u(x)=\sum_{i=1}^n c_i \psi_i(x)
\end{equation}
Using the same weight function as the basis function and integrating, we get:
\begin{equation}\label{eqn29}
	\int_{x_0}^{x_n} \psi_j (\alpha\frac{d^2u(x)}{dx^2}+\beta \frac{du(x)}{dx}+\gamma u(x) - f(x)) \,dx =0 
\end{equation}
On substituting equation \ref{eqn28}, using linear basis functions, and simplifying\cite{atluri2005methods}, equation \ref{eqn29} converts to the following set of linear equations.\\
\small 
\begin{equation*}
	\begin{bmatrix}
		(-\frac{\alpha}{h}+\frac{\alpha x_1}{h^2} -\frac{\beta}{2} -\frac{\gamma h}{3} -\frac{\alpha x_2}{h^2}) & (-\frac{\alpha x_1}{h^2}+\frac{\beta}{2}-\frac{\gamma h}{6}+\frac{\alpha x_2}{h^2}) & ... &0\\
		(-\frac{\alpha x_1}{h^2} -\frac{\beta}{2}+\frac{\gamma h}{6} +\frac{\alpha x_2}{h^2}) & 
		(\frac{\alpha x_1}{h^2} +\frac{2\gamma h}{3}-\frac{\alpha x_3}{h^2}) & (-\frac{\alpha x_2}{h^2}+\frac{\beta}{2}+\frac{\gamma h }{6}+\frac{\alpha x_3}{h^2})&0\\
		& \vdots & ... & \\
		0& (-\frac{\alpha x_{i-1}}{h^2}-\frac{\beta}{2} +\frac{\gamma h}{6}+\frac{\alpha x_i}{h^2}) & 
		(\frac{\alpha x_{i-1}}{h^2} +\frac{2\gamma h}{3} -\frac{\alpha x_{i+1}}{h^2}) & (-\frac{\alpha x_i}{h^2} +\frac{\beta}{2} +\frac{\gamma h}{6}+\frac{\alpha x_{i+1}}{h^2})\\
		0&0&(-\frac{\alpha x_{n-1}}{h^2} -\frac{\beta}{2}+\frac{\gamma h}{6} -\frac{\alpha x_n}{h^2}) & (\frac{\alpha x_{n-1}}{h^2} +\frac{\beta}{2}+\frac{\gamma h}{3}-\frac{\alpha x_{n}}{h^2})
	\end{bmatrix}
	\begin{bmatrix}
		c_1\\c_2\\c_3\\ \vdots\\c_n
	\end{bmatrix}		
\end{equation*}
\begin{equation}\label{30}
	=  \begin{bmatrix}
		1/h\int_{x_1}^{x_2}(x_2-x)f(x)\\ \vdots \\1/h\int_{x_{i-1}}^{x_i}(x-x_{i-1})f(x) + 1/h\int_{x_{i}}^{x_{i+1}}(x_{i+1}-x)f(x)\\ \vdots\\ 1/h\int_{x_{n-1}}^{x_n}(x-x_{n-1})f(x) - c\alpha
	\end{bmatrix}
\end{equation}
\normalsize

For evaluating the integrals, Quadrature methods are commonly used. Since the integrands are repeated, we propose using DNNI for substantial speedup. Once the primitive is approximated in the given domain, all integrals can be obtained instantaneously by just changing the limits. Let the anti-derivatives of $f(x)$ and $xf(x)$ using DNNI are $N_1(x)$ and $N_2(x)$. A typical integral on the right-hand side of equation \ref{30} can be calculated as
\begin{equation*}
	\int_{x_{i}}^{x_{i+1}}(x_{i+1}-x)f(x) = x_{i+1}\int_{x_{i}}^{x_{i+1}}f(x)\,dx - \int_{x_{i}}^{x_{i+1}}xf(x)\,dx
\end{equation*}
\begin{equation}
	= x_{i+1}(N_1(x_{i+1})-N_1(x_i)) - (N_2(x_{i+1})-N_2(x_i))
\end{equation}
The above modification will significantly reduce the computation time for a higher number of nodes.
\begin{prop}\label{th1}
	For differential equations with source terms, there is a finite value for the number of nodes 'n,' after which the DNNI-based Galerkin method is computationally less expensive than the naive Galerkin method.  
\end{prop}
\begin{proof}
  	Let the differential equation with a source term be:
	\begin{equation}\label{gde}
		F(\frac{\partial y}{\partial t},...\frac{\partial^2y}{\partial x^2},\frac{\partial y}{\partial x},.. y,t...x) = S(x)
	\end{equation}
Using simplifications similar to equation \ref{eqn29}, it can be converted to a system of linear equations like equation \ref{30}. The constructed equation will be similar to:
\begin{equation*}
	\begin{bmatrix}
		 Matrix\ depending \\ on\ the\ LHS \\ of\ equation\ \ref{gde}
	\end{bmatrix} \begin{bmatrix}
		c_1\\c_2\\c_3\\ \vdots\\c_n
	\end{bmatrix} =
\end{equation*}
\begin{equation}\label{matrix eqn}
	 \begin{bmatrix}
		1/h\int_{x_1}^{x_2}f_1(x,x^2...x_1,x_2...)S(x)+ some\ terms\\ \vdots \\1/h\int_{x_{i-1}}^{x_i}f_2(x,x^2...x_1,x_2...)S(x) + 1/h\int_{x_{i}}^{x_{i+1}}f_3(x,x^2...x_1,x_2...)S(x)\\ \vdots\\ 1/h\int_{x_{n-1}}^{x_n}f_{2n-2}(x,x^2...x_1,x_2...)S(x) + some\ terms
	\end{bmatrix}
\end{equation}

Though using a clever implementation of DNNI, even the LHS of equation \ref{matrix eqn} can be computed efficiently; we are ignoring this because it often forms patterns and does not require repeated integrations. For the naive Galerkin method, numerical techniques of definite integrations are commonly used. For the 'n' number of nodes, the number of integrations performed is $2n-2$. Assuming that the fastest numerical technique takes '$t$' time to compute a single integral. Thus, the minimum time required to compute the RHS of equation \ref{matrix eqn} is $t(2n-2)$. If the rest of the method takes time '$\tau$,' the total time required for the naive Galerkin method is:
\begin{equation}
	t_1 = t(2n-2)+\tau(n)
\end{equation} 
Using DNNI the RHS term can be computed using combinations of a finite number of anti-derivatives such as $\int xS(x)\,dx$, $\int x^2S(x)\,dx$...., $\int S(x)\,dx$. Let the number of such anti-derivatives be '$m$,' and the average time taken for a DNNI approximation is '$T$ .'Thus, the total time to compute all anti-derivatives is approximately '$mT$ .'Also, let the time for calculating the definite integrals by putting limits on the DNNI anti-derivatives is '$\epsilon$ .'So, 
\begin{equation*}
	Total\ time= time\ taken\ for\ \left( DNNI\ +\ applying\ limits\ +\ rest\ of\ the\ process\right) . 
\end{equation*}
\begin{equation}
	\implies t_2	=\ mT + (2n-2)m\epsilon+ \tau(n)
\end{equation}
Now, 
\begin{equation}
	\frac{t_1}{t_2}=\frac{t(2n-2)+\tau(n)}{mT + (2n-2)m\epsilon+\tau(n)}
\end{equation}
is the ratio of time taken by the naive Galerkin method and the DNNI-based Galerkin method. For,
\begin{equation*}
	t_1 > t_2
\end{equation*}
\begin{equation*}
	\implies t(2n-2)+\tau(n)>mT + (2n-2)m\epsilon+\tau(n)
\end{equation*}
\begin{equation}\label{n}
	\implies n > 1+\frac{mT}{2(t-m\epsilon)}
\end{equation}
The value of $m$ depends on the type of basis function used. Complex basis functions will lead to several different integrals, and the value of $m$ will increase.
Since $\epsilon$ is much less than $t$, if $m$ is not too large, we will always get a finite value '$n_{critical}$' of the number of nodes, after which the DNNI-based Galerkin method will compute faster.  
\end{proof}
In the following test case, the values observed are:
\begin{center}
\begin{tabular}{c c}
	\hline\hline
	Parameters & Value \\ \hline
	\hline 
	$T$ & 2.810464692115784  \\\hline
	$t$ & 0.00010534977912902833 \\\hline
	$\epsilon$ & 5.729198455810547e-07\\\hline
	$m$ & 2\\\hline
	\hline
\end{tabular}	
\end{center}
Using equation \ref{n}, the theoretical break-even value of 'n' is,
\begin{equation}\label{quad}
n_{critical} = 1+\frac{2\times2.810464692}{2(0.00010534-2\times5.7291984e-07)} =  26971
\end{equation}
\begin{figure}[h]
\centering
\includegraphics[width=12cm,height=8.5cm]{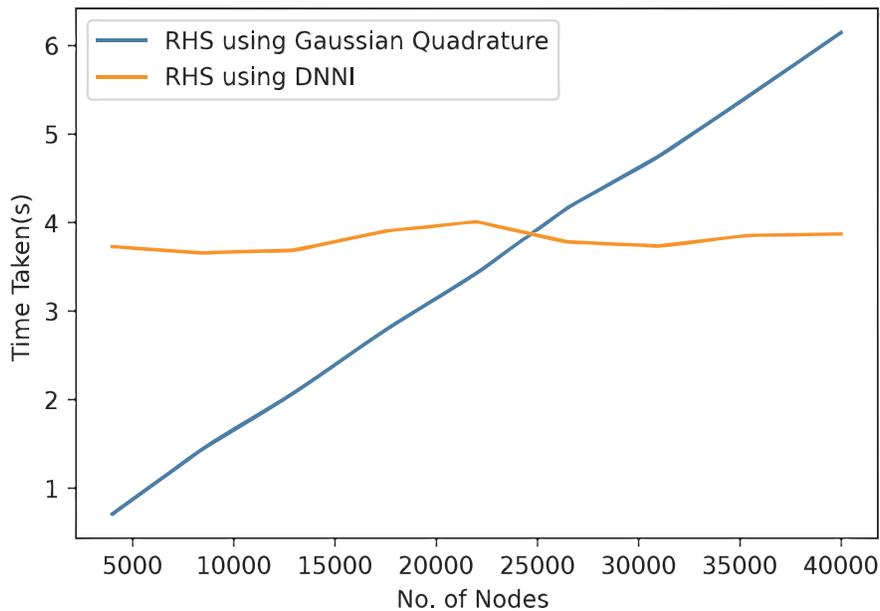}
\caption{Comparison of computation time for the RHS of equation \ref{30} using Gaussian Quadrature and DNNI method. The value of the critical number of nodes is of the same order as predicted in equation \ref{quad}.}
\label{fig:quad_dnni}
\end{figure}

\begin{table}
\begin{center}
	\begin{tabular}{c c c}
		\hline\hline
		Parameters & Case 14 & Case 15 \\ \hline
		\hline 
		$\alpha$ & 0 &0 \\\hline
		$\beta$ & 1&1 \\\hline
		$\gamma$ & 0&1\\\hline
		$u_o$ & 0&0\\\hline
		f(x) & cos(2x)&cos(2x)\\
		\hline \hline
	\end{tabular}	
\end{center}
\caption{The above parameters are used to solve two differential equations using the quadrature and DNNI-based Galerkin method. }
\end{table}

The speedup achieved is not very significant for simple test cases, but we conjecture that for more complex cases with computationally intensive repeated integrations, the DNNI-based method will outperform traditional quadrature-based Galerkin methods. 

\begin{figure}[h!]
\centering
\includegraphics[width=13cm,height=9.5cm]{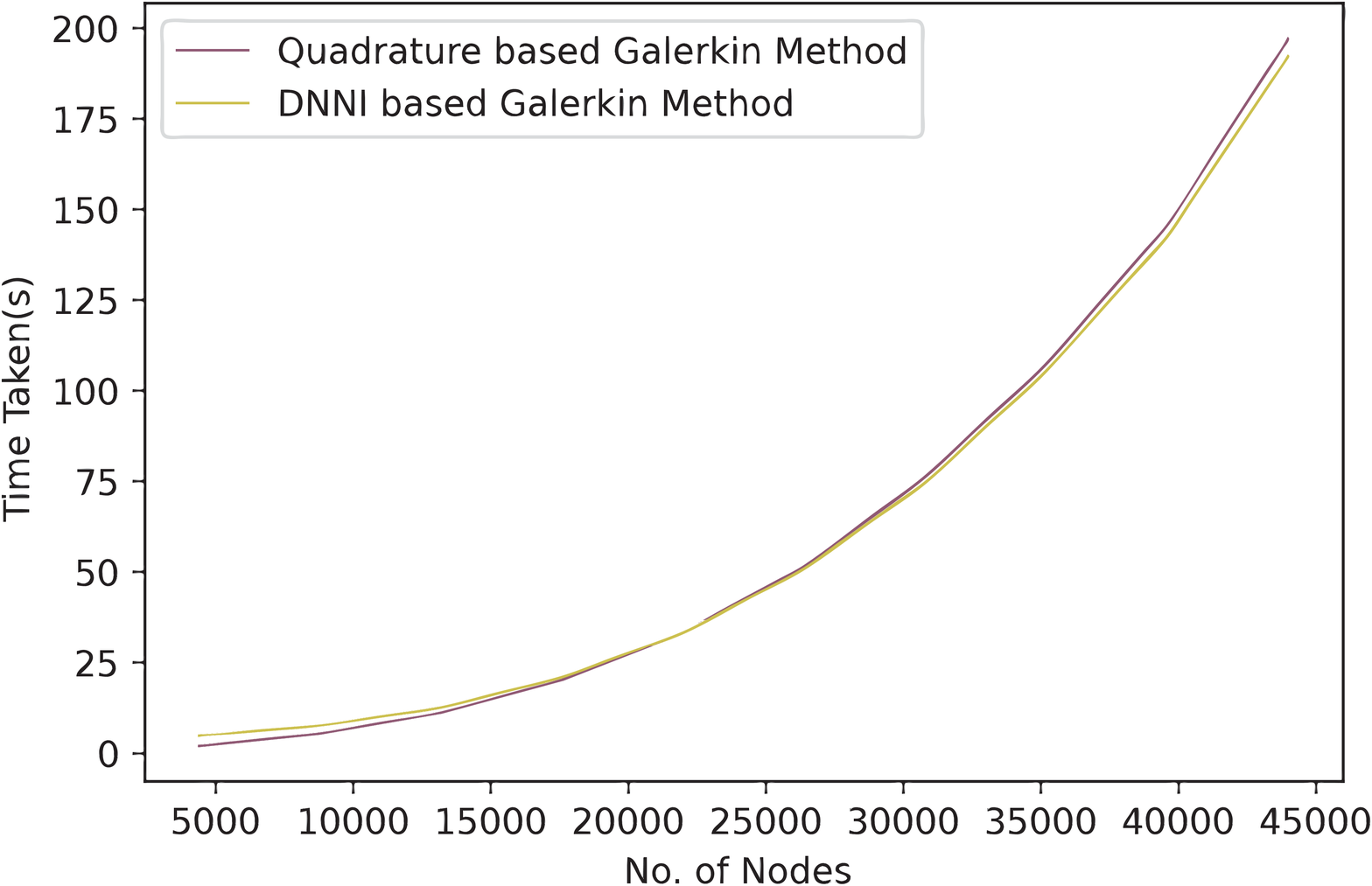}
\caption{Comparison of Quadrature based Galerkin method and DNNI based Galerkin method for Case 14. There is only a slight speedup because solving the set of linear equations is the most time-consuming step. However, the theoretical break-even point is of the same magnitude as the one obtained by computation.}
\label{fig: Galerkin}
\end{figure}
\begin{figure}[h!]
	\centering
	\includegraphics[width=13cm,height=9.5cm]{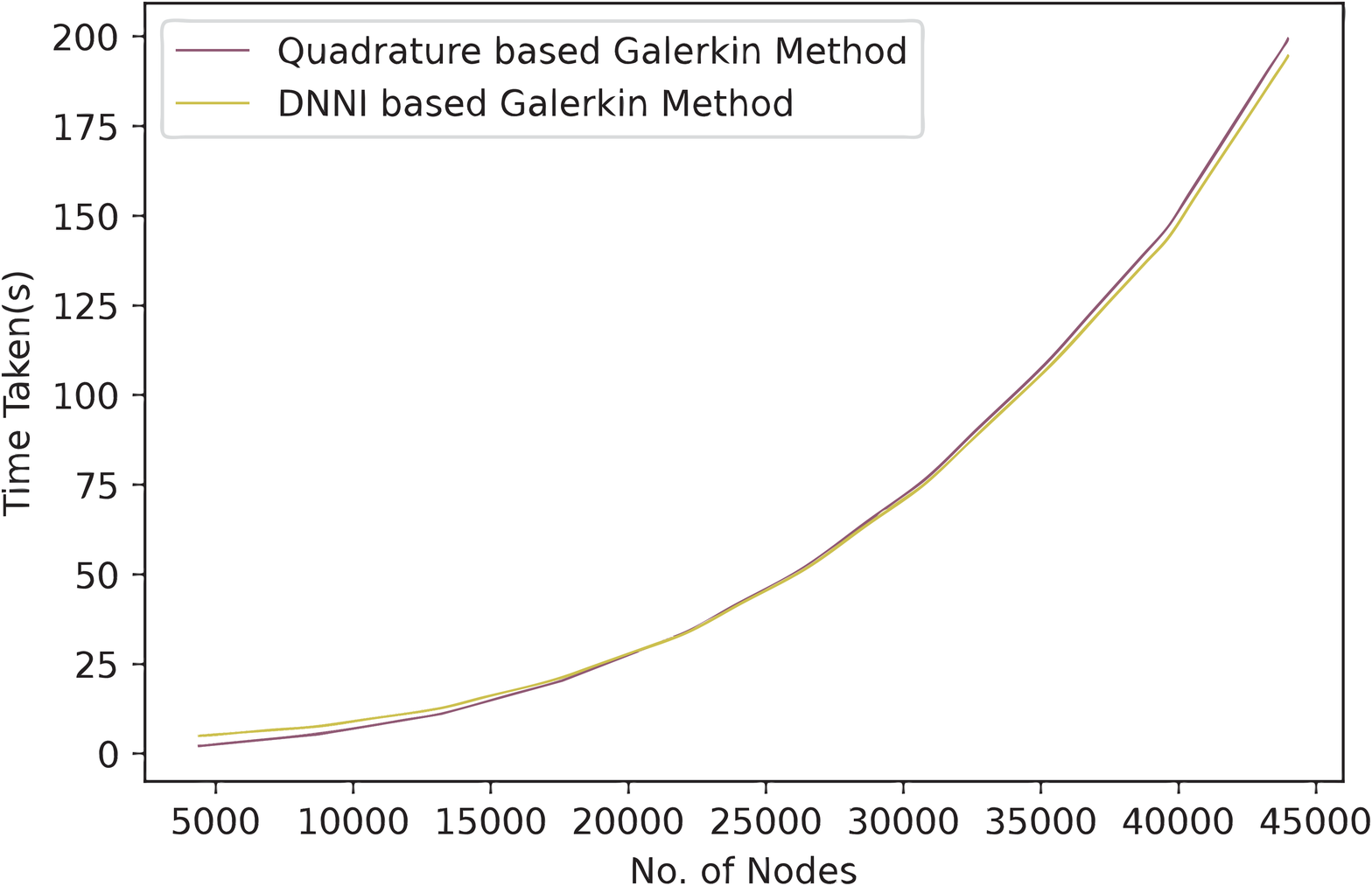}
	\caption{Comparison of Quadrature based Galerkin method and DNNI based Galerkin method for Case 15. }
	\label{fig: Galerkin2}
\end{figure}
%\clearpage
\newpage
\section{Conclusion}
In this paper, we propose an algorithm to represent the anti-derivative of a function using Deep Neural Networks. We have shown that DNNI effectively approximates very complex, non-elementary, and oscillatory integrals. We have also used DNNI to obtain the parameterized closed-form integrals, which can be later utilized to study the effects of various parameters. The closed form representations of the Fermi-Dirac and elliptic integrals were computed with significant accuracy.\\
Cumulative distribution functions were obtained using DNNI, eliminating the need for repeated numerical integrations to get the CDF tables once the anti-derivative is approximated. Also, the integrand need not be represented as a continuous function for DNNI. It can give a closed form anti-derivative on inputting even discrete values as integrand. Another advantage of DNNI is that it will get faster and more accurate with new optimization algorithms and Neural Network architecture developments.\\
The computational speedup using DNNI is shown for the Galerkin method, where repeated integrations are performed. DNNI can theoretically outperform quadrature-based Galerkin methods since it can instantly compute all the definite integrals after the anti-derivatives are obtained. Test cases were also computed to verify this proposition. The only downside of this method is that all the integration terms have to be reduced to a few anti-derivatives, which can be used repeatedly. This paper shows a clever approach to overcome this limitation for Galerkin methods with linear basis functions.\\
Further research on DNNI can include its application in the complex domain and multi-variable integrals. DNNI can also be applied to full-scale engineering problems where repeated or non-elementary integrals are required for obtaining substantial speedup. The application of DNNI on Galerkin methods for both temporal and spatial domains and using higher order polynomial basis functions is also an area of further research. DNNI can also be applied to compute the closed-form expressions of non-elementary integrals appearing in various areas of science. 

%Bibliography
\bibliographystyle{unsrt}  
\bibliography{arxiv}

\begin{thebibliography}{10}

\bibitem{muhammad2003double}
Mayinur Muhammad and Masatake Mori.
\newblock Double exponential formulas for numerical indefinite integration.
\newblock {\em Journal of Computational and Applied Mathematics},
  161(2):431--448, 2003.

\bibitem{hornik1989multilayer}
Kurt Hornik, Maxwell Stinchcombe, and Halbert White.
\newblock Multilayer feedforward networks are universal approximators.
\newblock {\em Neural networks}, 2(5):359--366, 1989.

\bibitem{zhe2006numerical}
Zeng Zhe-Zhao, Wang Yao-Nan, and Wen Hui.
\newblock Numerical integration based on a neural network algorithm.
\newblock {\em Computing in science \& engineering}, 8(4):42--48, 2006.

\bibitem{li2019dual}
Haibin Li, Yangtian Li, and Shangjie Li.
\newblock Dual neural network method for solving multiple definite integrals.
\newblock {\em Neural computation}, 31(1):208--232, 2019.

\bibitem{lloyd2020using}
Steffan Lloyd, Rishad~A Irani, and Mojtaba Ahmadi.
\newblock Using neural networks for fast numerical integration and
  optimization.
\newblock {\em IEEE Access}, 8:84519--84531, 2020.

\bibitem{baydin2018automatic}
Atilim~Gunes Baydin, Barak~A Pearlmutter, Alexey~Andreyevich Radul, and
  Jeffrey~Mark Siskind.
\newblock Automatic differentiation in machine learning: a survey.
\newblock {\em Journal of Marchine Learning Research}, 18:1--43, 2018.

\bibitem{kingma2014adam}
Diederik~P Kingma and Jimmy Ba.
\newblock Adam: A method for stochastic optimization.
\newblock {\em arXiv preprint arXiv:1412.6980}, 2014.

\bibitem{lample2019deep}
Guillaume Lample and Fran{\c{c}}ois Charton.
\newblock Deep learning for symbolic mathematics.
\newblock {\em arXiv preprint arXiv:1912.01412}, 2019.

\bibitem{bronstein1998symbolic}
Manuel Bronstein.
\newblock Symbolic integration tutorial.
\newblock Citeseer, 1998.

\bibitem{shivaram2016numerical}
KT~Shivaram and HT~Prakasha.
\newblock Numerical integration of highly oscillating functions using
  quadrature method.
\newblock {\em Global Journal of Pure and Applied Mathematics}, 3:2683--2690,
  2016.

\bibitem{gil2022complete}
Amparo Gil, Javier Segura, and Nico~M Temme.
\newblock Complete asymptotic expansions for the relativistic fermi-dirac
  integral.
\newblock {\em Applied Mathematics and Computation}, 412:126618, 2022.

\bibitem{sagar1991gaussian}
Robin~P Sagar.
\newblock A gaussian quadrature for the calculation of generalized fermi-dirac
  integrals.
\newblock {\em Computer physics communications}, 66(2-3):271--275, 1991.

\bibitem{temme1990uniform}
Nico~M Temme and AB~Olde Daalhuis.
\newblock Uniform asymptotic approximation of fermi—dirac integrals.
\newblock {\em Journal of Computational and Applied Mathematics},
  31(3):383--387, 1990.

\bibitem{bhagat2003evaluation}
Vikram Bhagat, Ranjan Bhattacharya, and Dhiranjan Roy.
\newblock On the evaluation of generalized bose--einstein and fermi--dirac
  integrals.
\newblock {\em Computer physics communications}, 155(1):7--20, 2003.

\bibitem{mohankumar2016very}
N~Mohankumar and A~Natarajan.
\newblock On the very accurate numerical evaluation of the generalized
  fermi--dirac integrals.
\newblock {\em Computer Physics Communications}, 207:193--201, 2016.

\bibitem{giraldo2020introduction}
Francis~X Giraldo.
\newblock {\em An Introduction to Element-Based Galerkin Methods on
  Tensor-Product Bases: Analysis, Algorithms, and Applications}, volume~24.
\newblock Springer Nature, 2020.

\bibitem{atluri2005methods}
Satya~N Atluri.
\newblock {\em Methods of computer modeling in engineering \& the sciences},
  volume~1.
\newblock 2005.

\end{thebibliography}

\end{document}